\makeatletter
\def\cl@chapter{}
\makeatother

\documentclass[smallextended,nospthms]{svjour3}       % onecolumn (second format)
\smartqed  % flush right qed marks, e.g. at end of proof
\usepackage{graphicx}

\usepackage{mathptmx}      % use Times fonts if available on your TeX system
\usepackage[square,sort,comma,numbers]{natbib}

\usepackage[all]{xy}
\usepackage{subfigure}

\usepackage{cleveref}
\usepackage{graphicx}
\usepackage[ruled]{algorithm2e} 
\usepackage[font=footnotesize]{caption}

\usepackage{amsmath}
\DeclareMathOperator{\arccosh}{acosh}
\usepackage{amsthm} 
\usepackage{mathtools} 
\usepackage{hyperref}
\usepackage{etoolbox}

\newtheorem{defn}{Definition}[section]

\newtheorem{rem}{Remark}
\newcommand{\x}{\textbf{x}}
\newcommand{\z}{\textbf{z}}
\newcommand{\xtilde}{\tilde{\x}}
\newcommand{\xhat}{\hat{\x}}
\newcommand{\ba}{\textbf{a}}
\newcommand{\bb}{\textbf{b}}
\newcommand{\bc}{\textbf{c}}
\newcommand{\bv}{\textbf{v}}

\newcommand{\VN}{{\mathchoice{}{}{\scriptscriptstyle}{}V}}
\newcommand{\HN}{{\mathchoice{}{}{\scriptscriptstyle}{} H}}
\newcommand{\uH}{u_{\HN}^{}}
\newcommand{\uV}{u_{\VN}^{}}
\newcommand{\kH}{k_{\HN}^{}}
\newcommand{\kV}{k_{\VN}^{}}
\newcommand{\tauH}{\tau_{\HN}^{}}
\newcommand{\tauV}{\tau_{\VN}^{}}

\newtoggle{ForJournal}
\toggletrue{ForJournal}

\journalname{Journal of Scientific Computing}

\begin{document}

\title{Corner cases, singularities, and dynamic factoring\thanks{
The second author's work is supported in part by the National Science Foundation grant DMS-1738010
and the Simons Foundation Fellowship.}
}

\titlerunning{Dynamic Factoring}        % if too long for running head

\author{Dongping Qi         \and
        Alexander Vladimirsky %etc.
}


\institute{D. Qi \at
              Center for Applied Mathematics, Cornell University, Ithaca, NY 14853\\
              \email{dq48@cornell.edu}           %  \\
           \and
           A. Vladimirsky \at
              Dept. of Mathematics and Center for Applied Mathematics, Cornell University, Ithaca, NY 14853\\
              \email{vladimirsky@cornell.edu}
}

\date{Received: date / Accepted: date}

\maketitle

%-----------------------------------------------------------------

% !TEX root = local_factoring_for_Journal.tex

\begin{abstract}
\noindent
In Eikonal equations, rarefaction is a common phenomenon known to degrade the rate of convergence of numerical methods.  The ``factoring'' approach alleviates this difficulty by deriving a PDE for a new (locally smooth) variable while capturing the rarefaction-related singularity in a known (non-smooth) ``factor''.  Previously this technique was successfully used to address rarefaction fans arising at point sources.  In this paper we show how similar ideas can be used to factor the 2D rarefactions arising due to nonsmoothness of  domain boundaries or discontinuities in PDE coefficients.  Locations and orientations of such rarefaction fans are not known in advance and we construct a ``just-in-time factoring'' method that identifies them dynamically.  The resulting algorithm is a generalization of the Fast Marching Method originally introduced for the regular (unfactored) Eikonal equations.  We show that our approach restores the first-order convergence and illustrate it using a range of maze navigation examples with non-permeable and ``slowly permeable'' obstacles.
\iftoggle{ForJournal}{
}{}
\end{abstract}

%-----------------------------------------------------------------

% !TEX root = local_factoring_for_Journal.tex

\section{Introduction}
\label{s:intro}
\setlength{\parindent}{2em}

\begin{figure}[!htb]
\centering
\subfigure[]{
\label{f:hessian_v}
\includegraphics[width=0.46\textwidth,clip]{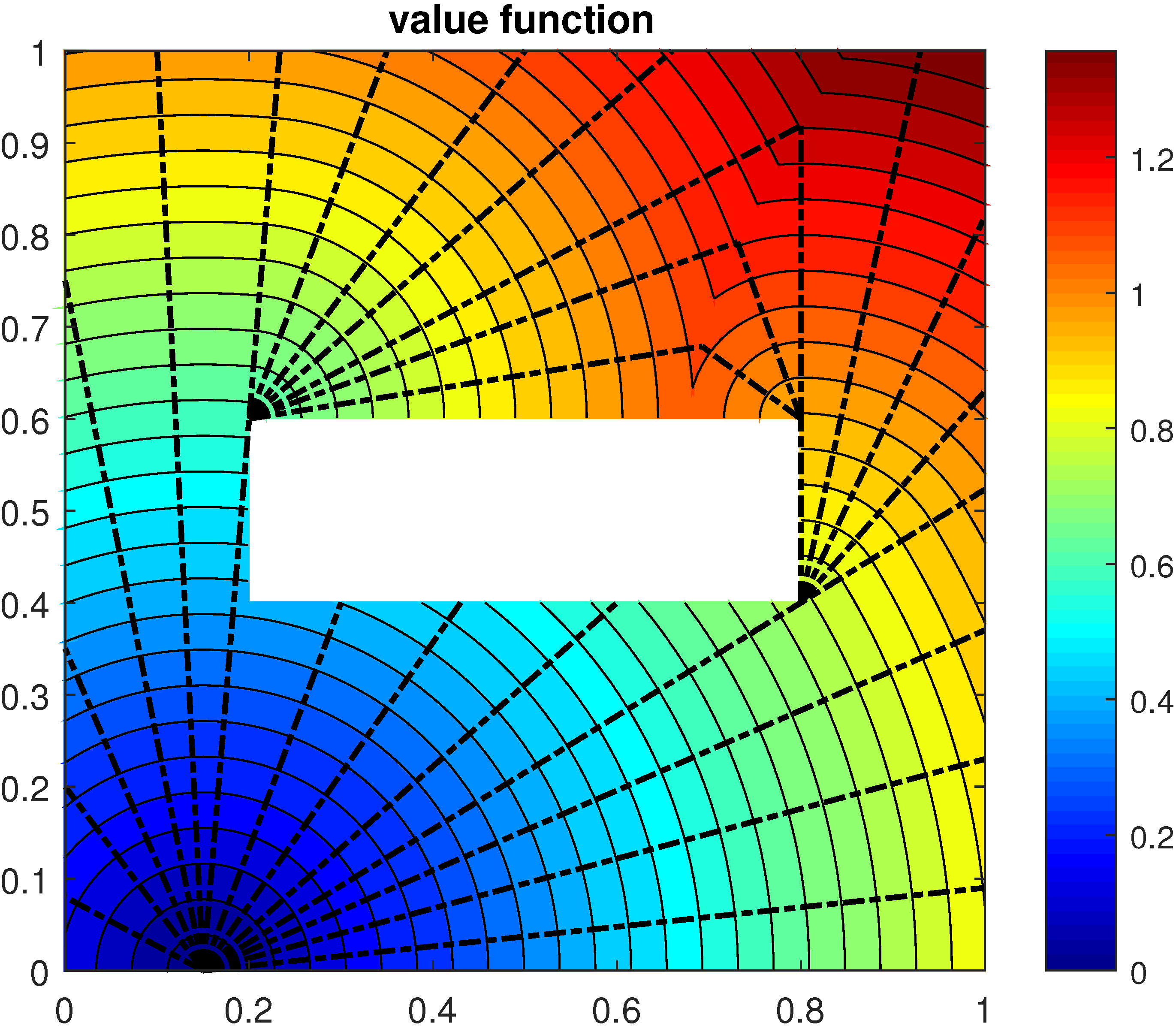}}
\subfigure[]{
\label{f:hessian_h}
\includegraphics[width=0.45\textwidth,clip]{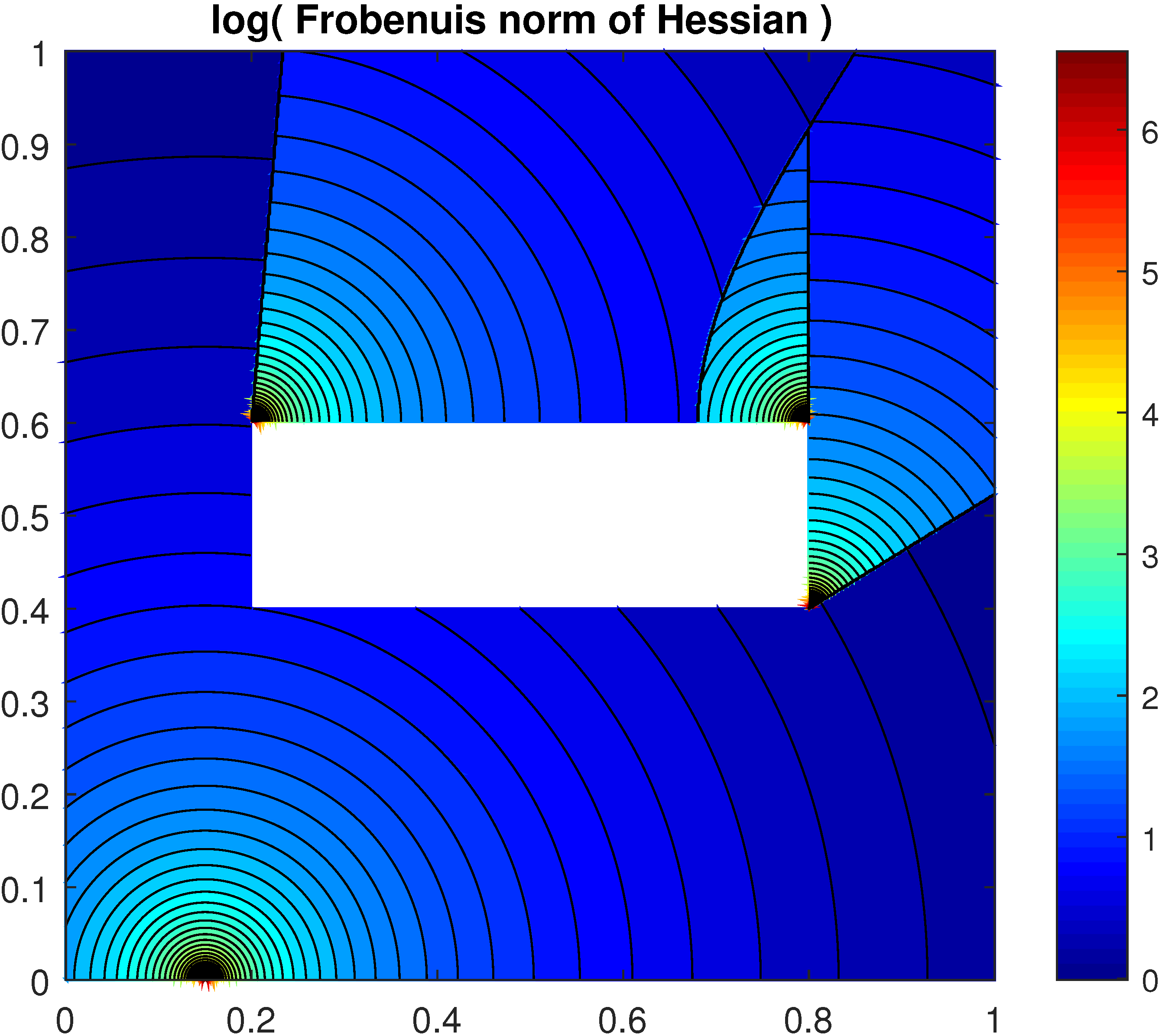}}
\caption{A simple example with $F=1$ and the exit set $Q$ consisting of a single ``point source'' $\x_0 = (0.15, 0).$
(Left) Level sets of $u(\x)$, the distance to $\x_0$  on the domain with an obstacle. 
Dashed lines show the optimal/shortest paths, some of which follow a part of the obstacle boundary.
The ``shockline'' (where $\nabla u$ is undefined) consists of all points starting from which the shortest path is not unique.
The point source and 3 out of 4 obstacle corners generate ``rarefaction fans'' of characteristics.
(Right) The level curves of the log-Frobenius-norm of Hessian, illustrating the blow up of the second derivatives of $u$ due to rarefaction fans.
}
\label{f:hessian}
\end{figure}

The Eikonal equation arises in many application domains including the geometric optics, optimal control theory, and robotic navigation. It is a first-order non-linear partial differential equation of the form
\begin{equation}
\label{eq:Eikonal}
\left\{
\begin{aligned}
 \left| \nabla u(\x) \right| F(\x) &= 1, &\\  
  u(\x) &= 0, & \x & \in Q.\\
\end{aligned}
\right.
\end{equation}
In control-theoretic context, the function $u(\x)$ can be interpreted as the minimum time needed to reach the exit set $Q$ 
starting from $\x$, with 
$F$ specifying the speed of motion.  
The characteristics of this PDE coincide with the gradient lines of $u$ and yield min-time-optimal trajectories to  $Q$.
This equation typically does not have a classical solution, making it necessary to select a  {\em weak} Lipschitz-continuous solution that is physically meaningful.  The theory of {\em viscosity solutions} \cite{crandall1983viscosity} accomplishes this and guarantees the existence and uniqueness for a very broad set of problems.  Viscosity solutions exhibit both shocklines (where characteristics run into each other) and rarefaction fans (where many characteristics emanate from the same point).  Typically, the former receive most of the attention in numerical treatment since they lead to a discontinuity in $\nabla u$.  But in this paper we focus on the latter, which result in a blow up in second derivatives of $u$ on parts of the domain where $\nabla u$ remains continuous and bounded.  (See Figure \ref{f:hessian}.)

Since the local truncation error of numerical discretizations is proportional to higher derivatives, it is not surprising that rarefaction fans degrade the rate of convergence of standard numerical methods.  
Special ``factoring techniques'' have been introduced to alleviate this problem for rarefaction fans arising from point source boundary conditions
\cite{fomel2009, luo2011factored, luo2012fast, luo2014high, noble2014, treister2016fast}.  
The idea is to represent $u$ as a product \cite{fomel2009} or a sum \cite{luo2012fast} of two functions: one capturing the right type of singularity at the point source and another with bounded derivatives except at shocklines.  The former is known based on $F$ restricted to $Q$; the latter is a priori unknown and recovered by solving a modified PDE by a version of one of the efficient methods (e.g., \cite{tstitsiklis1996, sethian1996fast} or \cite{zhao2005fast}) originally developed for an ``un-factored'' Eikonal equation \eqref{eq:Eikonal}.  We review this prior work in section \ref{s:point_sources} and show that, despite its better rate of convergence, factoring can also have detrimental effects on parts of 
the domain  far from the point source.
We also consider a  {\em localized} version of factoring, which often improves the accuracy and is more suitable for problems with multiple point sources.

However, point sources are not the only origin of rarefaction fans.  They can also arise due to non-smoothness of the boundary; e.g., at some of the ``obstacle corners'' in 2D maze navigation problems (see Figure \ref{f:hessian}).  
The degradation of accuracy leads to numerical artifacts in computed trajectories passing near such corners. 
Unlike in the point source case, these rarefaction fans are not radially symmetric; moreover, their locations and geometry have to be determined dynamically. 
We handle this by developing  a ``just-in-time localized factoring'' method and verifying its rate of convergence numerically (section \ref{s:bad_corners}).  Even more complicated fans can arise at corners  of ``slowly permeable obstacles'' (i.e., problems with piecewise-continuous speed 
function $F$).  The extension of our dynamic factoring to this case is covered in section \ref{ss:discont}.

Throughout the paper, we present our approach using a Cartesian grid discretization of the Eikonal equation with grid aligned obstacles (and discontinuities of $F$).  However, the ideas presented here have a broader applicability.
Arbitrary polygonal obstacles can be treated in exactly the same way if the discretization is posed on an obstacle-fitted triangulated mesh.  Dynamic factoring for more general Hamilton-Jacobi-Bellman PDEs  would work very similarly, though the factoring will need to account for the anisotropy in rarefaction fans.  We conclude by discussing these and other future extensions as well as the limitations of our approach in section \ref{s:conclude}.

%-----------------------------------------------------------------

% !TEX root = local_factoring_for_Journal.tex

\section{Point-sources and factoring}
\label{s:point_sources}

We begin by examining the classical Eikonal equation \eqref{eq:Eikonal} in 2D with a ``single point-source'' boundary condition:
$Q = \{ \x_0 \}$ and $u(\x_0) = 0$.
Throughout this paper, we will assume that the controlled process is restricted to a closed set $\Omega \subset R^2.$
I.e., $u(\x)$ is the minimal-time from $\x \in \Omega$ to $Q \subset \Omega$ without leaving $\Omega$ though possibly 
traveling along parts of $\partial \Omega$; see the trajectories traveling along the obstacle boundary in Figure \ref{f:hessian_v}.
This makes $u$ an {\em $\Omega$-constrained viscosity solution} of \eqref{eq:Eikonal} (see Chapter 4.5 in \cite{bardicapuzzo}),
but for the purposes of numerical implementation we can simply impose the boundary condition $u=+\infty$ on $R^2 \backslash \Omega$.

A common approach for discretizing \eqref{eq:Eikonal} on a uniform Cartesian grid is to use upwind finite differences:
\begin{equation}
\label{eq:scheme}
\max\{D_{i,j}^{-x}u, -D_{i,j}^{+x}u, 0\}^{2} + \max\{D_{i,j}^{-y}u, -D_{i,j}^{+y}u, 0\}^{2} = \frac{1}{F_{i,j}^2},
\end{equation}
using the standard finite difference notation
\begin{equation}
\label{eq:operator1}
\begin{split}
D_{i,j}^{-x}u = \frac{u_{i,j} - u_{i-1,j}}{h}, \qquad D_{i,j}^{+x}u= \frac{u_{i+1,j} - u_{i,j}}{h}, \\
D_{i,j}^{-y}u = \frac{u_{i,j} - u_{i,j-1}}{h}, \qquad D_{i,j}^{+y}u= \frac{u_{i,j+1} - u_{i,j}}{h}.
\end{split}
\end{equation}
The discretized system of equations \eqref{eq:scheme} is coupled and non-linear.
We postpone the discussion of fast algorithms used to solve it until section \ref{ss:fastEik}
and for now focus on the rate of convergence of the numerical solution to the viscosity solution of PDE \eqref{eq:Eikonal}.
Since this discretization is globally first-order accurate, the local truncation error is proportional to the second derivatives of $u$, which blow up as we approach $\x_0$.  Because of this blow up, even if we assume that local truncation errors accumulate linearly, the global error would decreases as $O(h \log\frac{1}{h})$ instead of the expected $O(h).$  

\begin{figure}[!htb]
\centering
\subfigure[level sets of $u(\x)$]{
\label{f:nonobst_1_v}
\iftoggle{ForJournal}{\includegraphics[width=0.47\textwidth,clip]{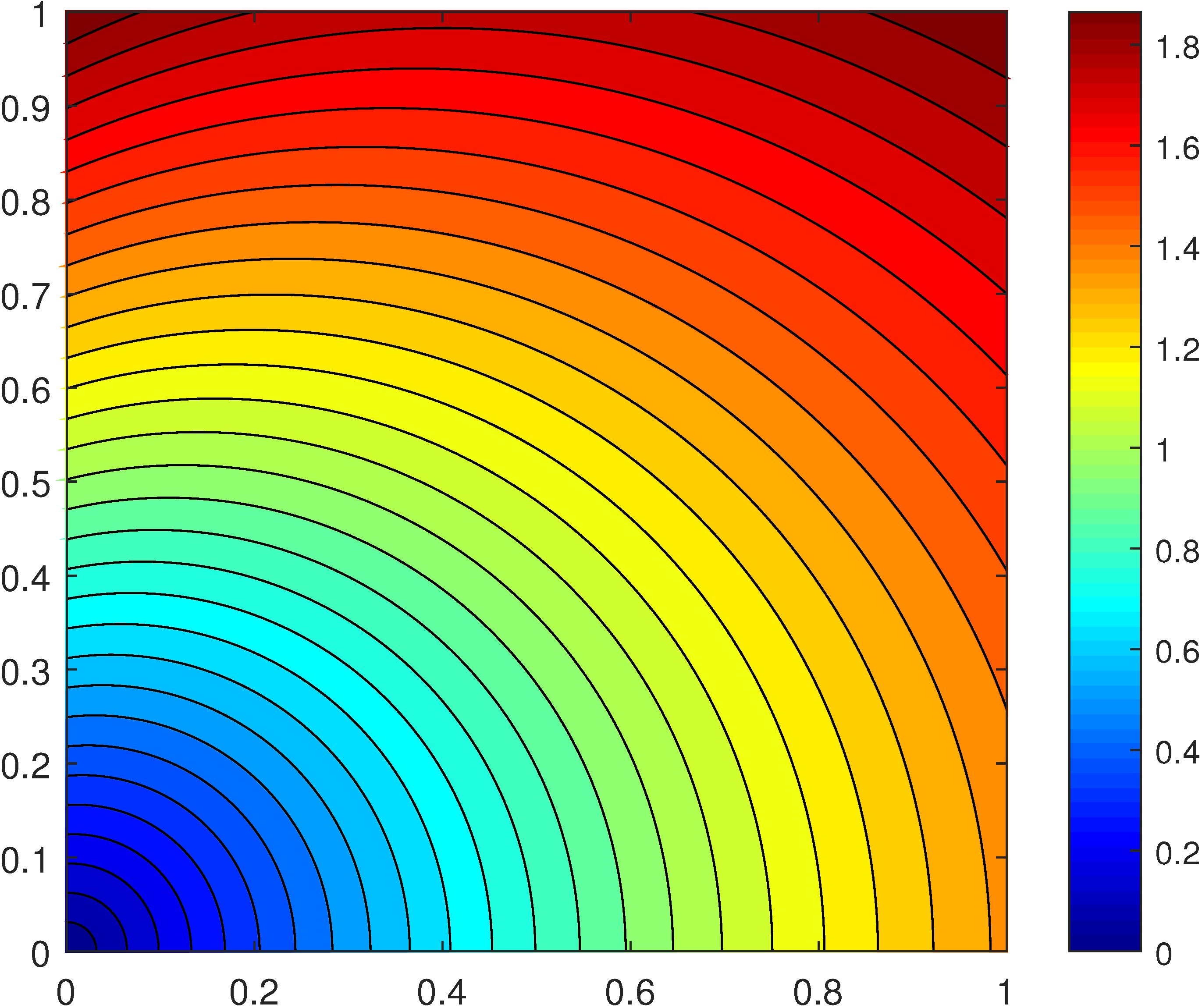}}{\includegraphics[width=0.38\textwidth,clip]{nonobst_1_v}}}
\subfigure[$L^{\infty}$ error]{
\label{f:nonobst_1_c0}
\hspace*{5mm}
\includegraphics[scale=.45,clip]{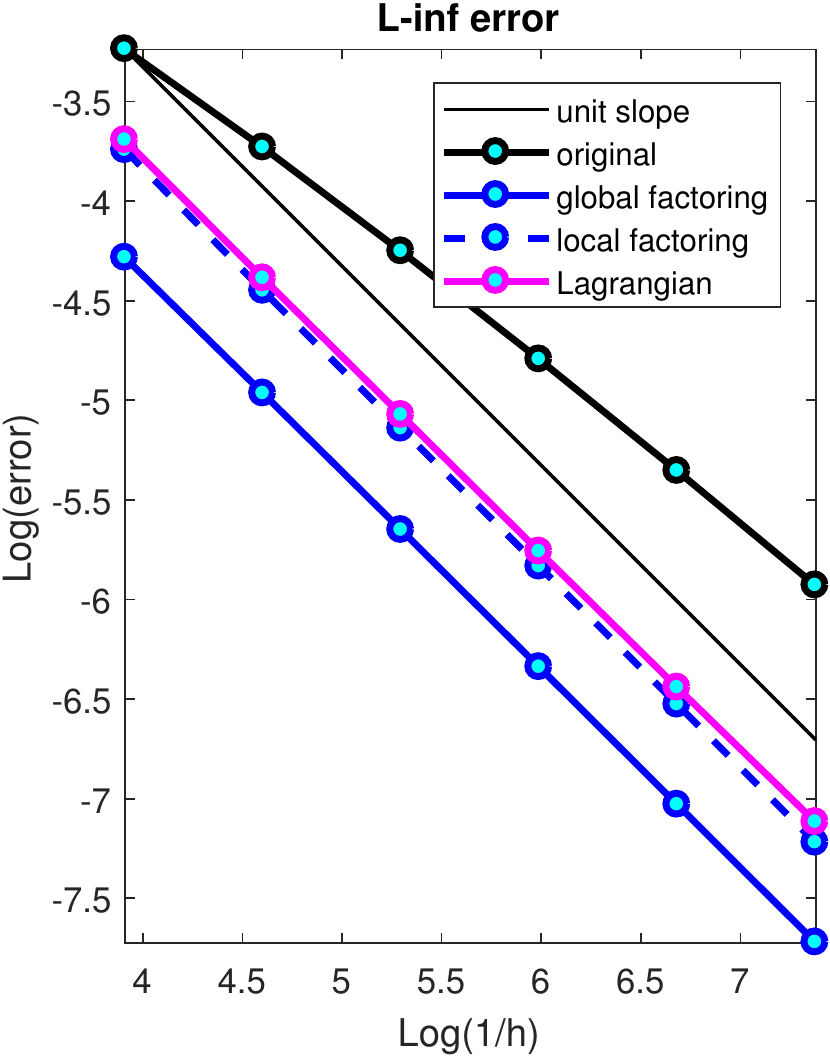}}
\caption{(Left) Level curves of the value function %$u(\x)$ 
computed from the formula \eqref{eq:single_src} with 
$\x_0=(0,0), \, s_0=2, \, \bv=(0.5,0).$
(Right) The corresponding convergence plot (based on the $L^{\infty}$ error) for several discretization approaches, with the thin black line of slope $(-1)$ included to aid the visual comparison.  For this example, the convergence plot based on the $L^1$ error is very similar and thus omitted.}
\label{f:nonobst_1}
\end{figure}

To illustrate the convergence rate, we will consider a simple example with a known analytic solution \cite{luo2012fast} on a square domain\footnote{
This formula for $u(\x)$ is derived on the unbounded domain $\Omega_{\infty} = \{ \x \in R^2 \mid F(\x) > 0 \}$ 
but remains valid 
on $\Omega = [0,1]\times[0,1]$ 
as long as $\Omega_{\infty}$-optimal trajectories from every $\x \in \Omega$ to $\x_0$ stay entirely inside $\Omega,$ which is the case for all examples considered in this section.
The linearity of $F(\x)$ can be used to show that all optimal paths are circular arcs, 
%TODO: ref?
whose radii are monotone decreasing in $|\bv|$. (When $\bv=0,$ these radii are infinite; i.e., all optimal paths are straight lines and 
$u(\x) = s_0 |\x-\x_0|.$)
}:
\begin{equation}
\label{eq:single_src}
F(\x) = \frac{1}{s_0} + \bv\cdot(\x-\x_0) \quad \Longrightarrow \quad
u(\x) = \frac{1}{|\bv|}\arccosh\left(1+\frac{1}{2}s_{0}{|\bv|}^2 \frac{\left| \x-\x_0 \right|^{2}}{F(\x)}\right).
\end{equation}
Figure \ref{f:nonobst_1} shows the solution and the convergence plots for the parameter values $\x_0=(0,0), \, s_0=2, \, \bv=(0.5,0).$
The thick black line in Figure \ref{f:nonobst_1_c0}
corresponds to solving  \eqref{eq:scheme} on $\Omega$ and clearly shows the sublinear convergence. 
(Throughout this paper, all logarithms of errors and grid sizes are reported and plotted in base $e.$) 

\begin{figure}[!htb]
\centering
\subfigure[level sets]{
\label{f:nonobst_2_v}
\includegraphics[width=0.38\textwidth,clip]{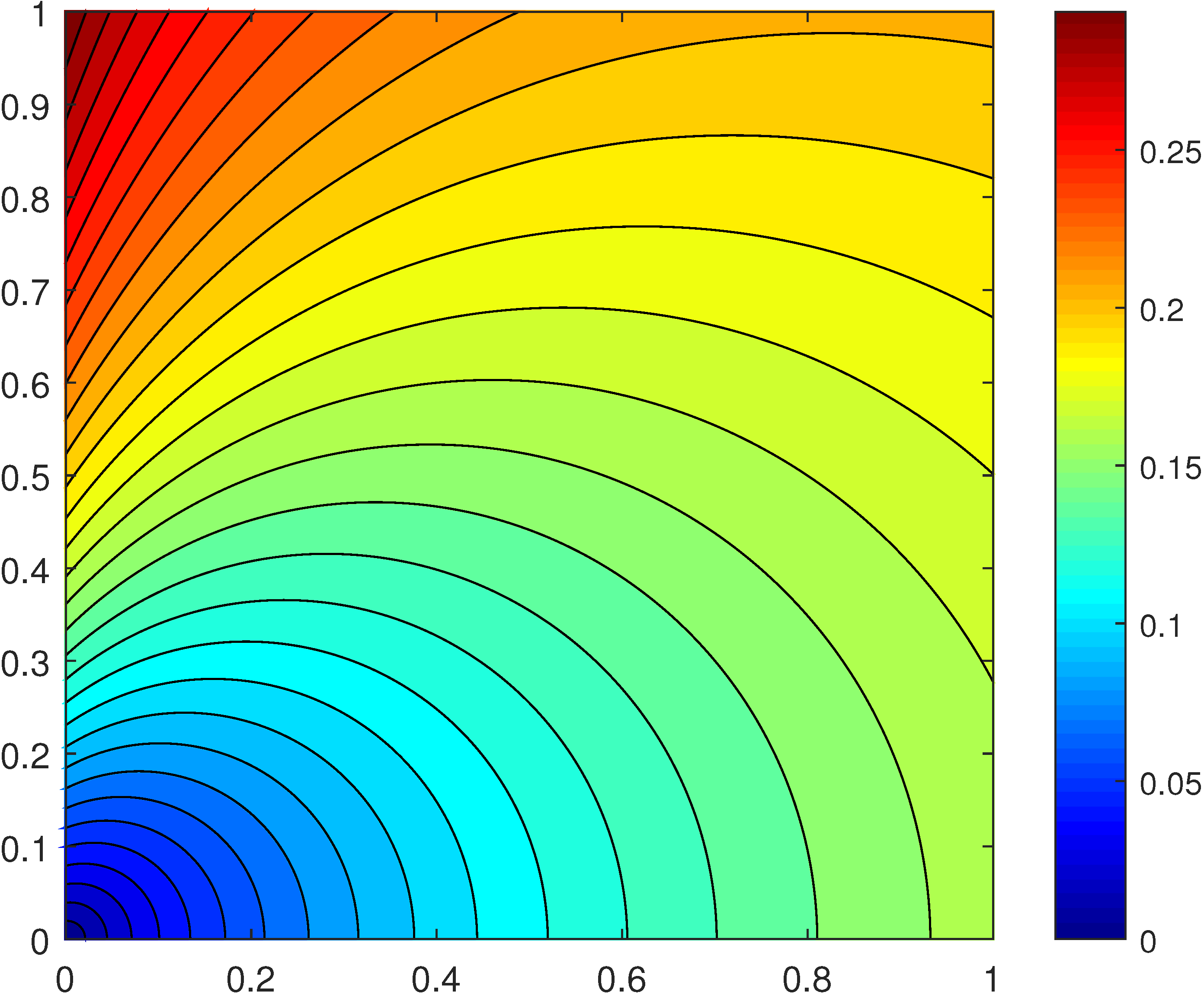}}
\subfigure[$L^{\infty}$ error]{
\label{f:nonobst_2_conv_Linf}
\iftoggle{ForJournal}{\includegraphics[scale=.36,clip]{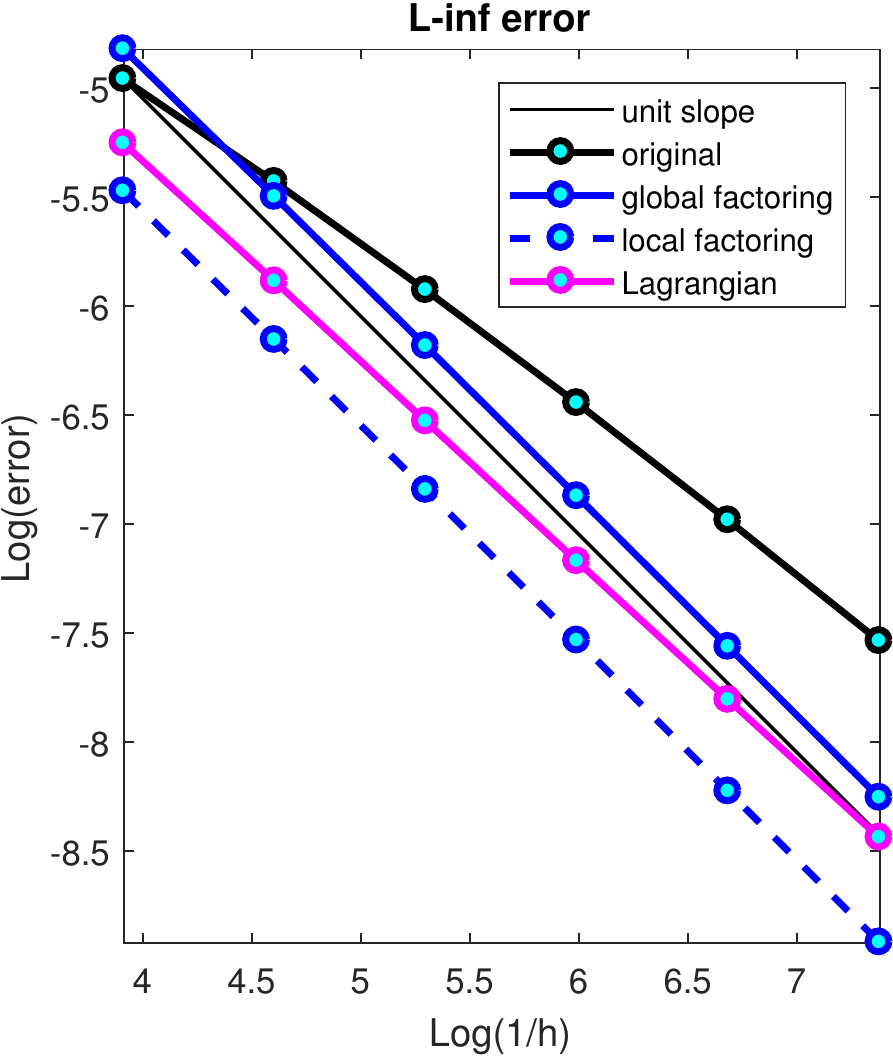}}{\includegraphics[scale=.45,clip]{nonobst_2_c0}}}
\subfigure[$L^1$ error]{
\label{f:nonobst_2_conv_L1}
\iftoggle{ForJournal}{\includegraphics[scale=.36,clip]{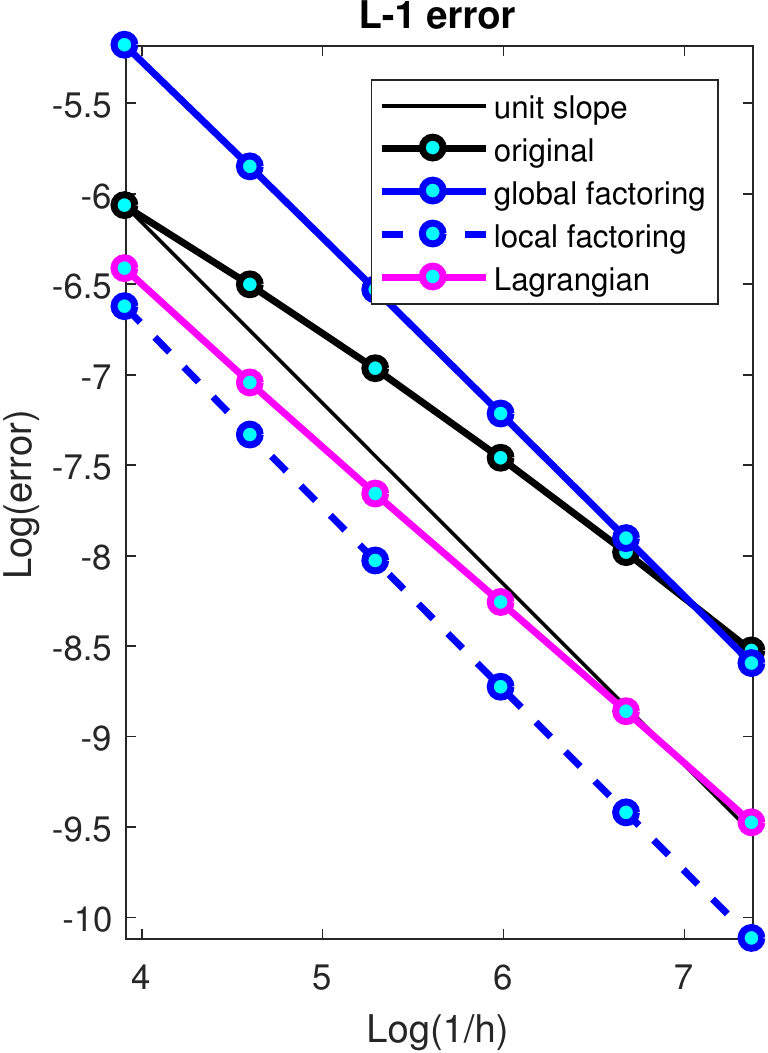}}{\includegraphics[scale=.45,clip]{nonobst_2_c1}}}
\caption{
(Left) Level curves of the value function $u(\x)$ computed from the formula \eqref{eq:single_src} with 
$\x_0=(0,0), \, s_0=0.5, \, \bv=(12,0).$
(Center \& Right) The corresponding convergence plot for several discretization approaches (based on the $L^{\infty}$ \& $L^1$ errors respectively) .}
\label{f:nonobst_2}
\end{figure}

\begin{figure}[!htb]
\centering
\subfigure[level sets]{
\label{f:nonobst_3_v}
\includegraphics[width=0.38\textwidth,clip]{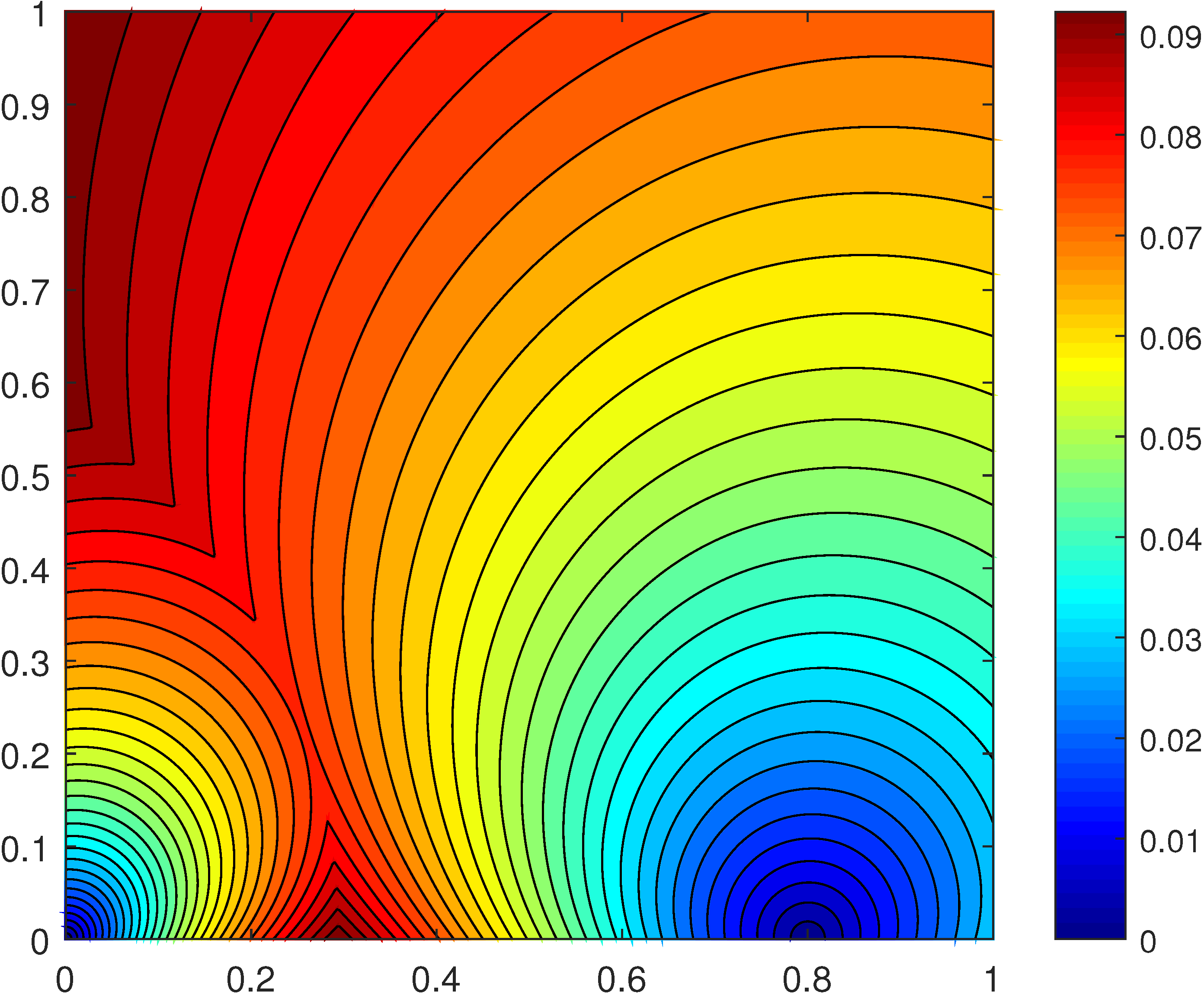}}
\subfigure[$L^{\infty}$ error]{
\label{f:nonobst_3_conv_Linf}
\iftoggle{ForJournal}{\includegraphics[clip,scale=.36]{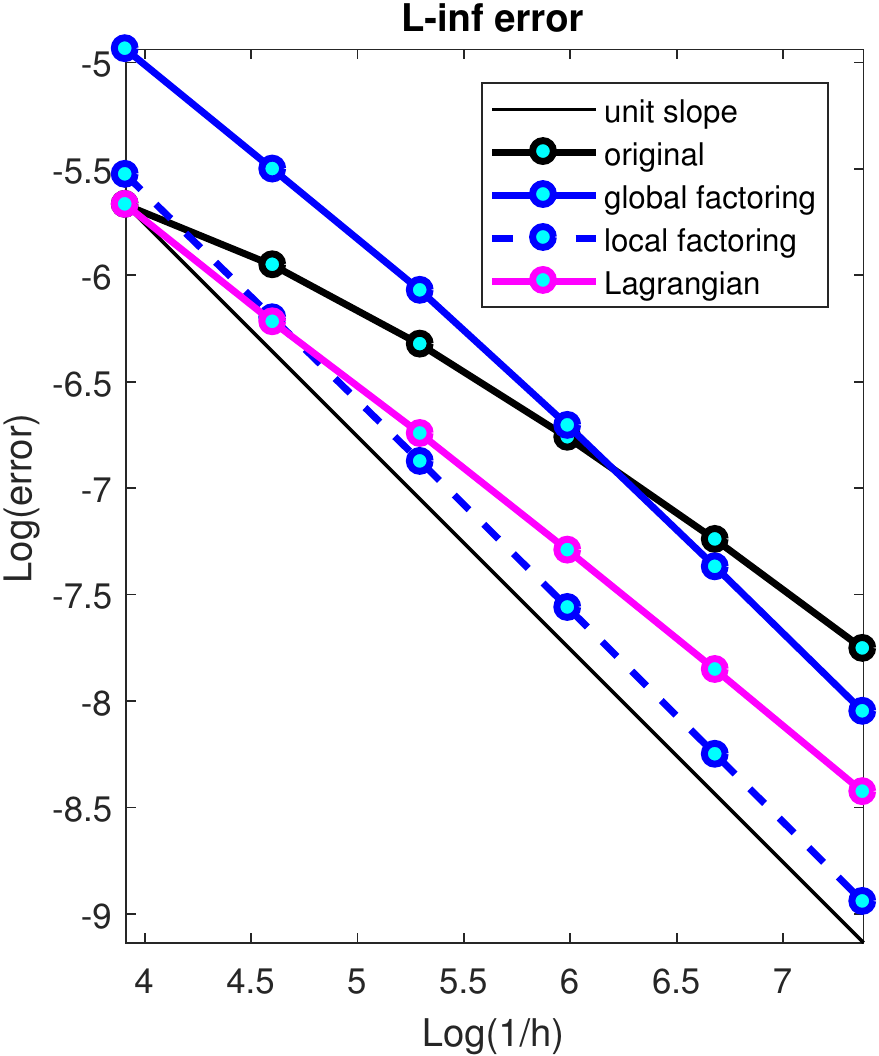}}{\includegraphics[clip,scale=.45]{nonobst_3_c0}}}
\subfigure[$L^1$ error]{
\label{f:nonobst_3_conv_L1}
\iftoggle{ForJournal}{\includegraphics[width=0.28\textwidth,clip]{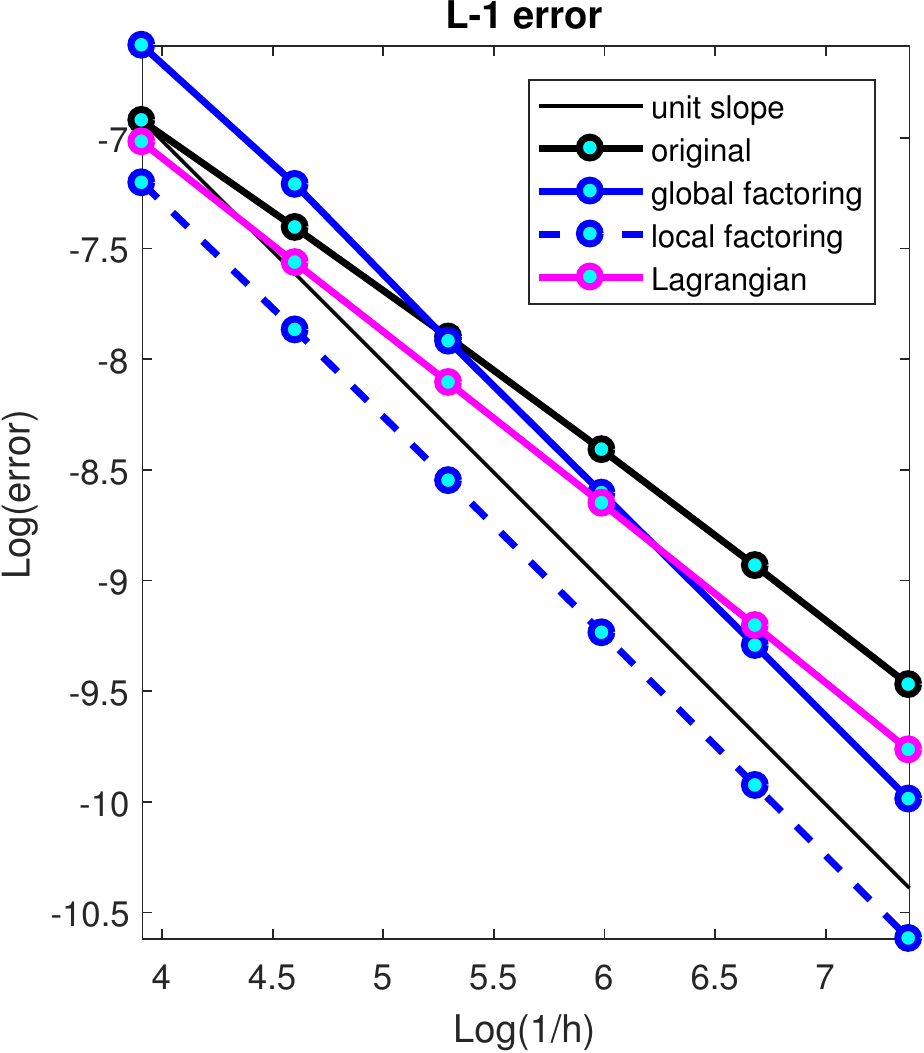}}{\includegraphics[width=0.28\textwidth,clip]{nonobst_3_c1}}}
\caption{
(Left) Level curves of the value function $u(\x)$ for $s_0 = 0.5, \, \bv = (5,20)$ and two point sources at $\x_0 = (0,0)$ and $\x_1=(0.8,0).$
(Center \& Right) The corresponding convergence plot for several discretization approaches (based on the $L^{\infty}$ \& $L^1$ errors respectively) .}
\label{f:nonobst_3}
\end{figure}

To alleviate this issue, one could simply ``enlarge the exit set'' by choosing some ($h$-independent) constant radius $r>0,$ initializing $u=0$ on the disk $B = B_r(\x_0) = \{ \x \in \Omega \mid \left| \x-\x_0 \right| \leq r \},$ and solving \eqref{eq:scheme} on the rest of the grid.
This avoids the rarefaction fan (since only one characteristic stems from each point on $\partial B$), but introduces a $O(r)$ difference compared to the solution of the original point-source-based problem.   One can also use a better approximation of $u$ on $B$; e.g., 
using $T(\x) = \frac{\left| \x-\x_0 \right|}{F(\x_0)}$, we already reduce this additional error to $O(r^2)$.  The latter approach is based on assuming that $F$ (rather than $u$) is constant on $B$, in which case the characteristics are straight lines.   Of course, one could also take a truly Lagrangian approach,
employing ray-tracing
to compute a highly accurate value of $u$ at all gridpoints on $B$, but this becomes increasingly expensive as $h \to 0$.  Instead, we evaluate the feasibility of an ``approximate Lagrangian'' initialization technique, where the characteristics are still assumed to be straight  on $B$, but $u(x)$ is approximated more accurately by integrating $1/F$ on the line segment from $\x_0$ to $\x$.
In all the figures of this section, we slightly abuse the notation and call this approach Lagrangian, reporting the results for $r=0.1$ in all the figures of this section.  
Figures \ref{f:nonobst_2} and \ref{f:nonobst_3} show that, for general $F$, the error due to this ``characteristics are straight'' assumption prevents the overall first order convergence on the entire grid.

{\em A factored Eikonal equation} was proposed in \cite{fomel2009} as a method for dealing with point-source rarefaction fans without introducing any special approximations on $B$ and recovering the first order of accuracy on the entire domain.
The main idea is to split the original value function $u(\x)$ into two functions: one of them ($T$, defined above) encodes the right type of singularity at the point source, while the other ($\tau(\x)$, our new unknown) will be smooth near $\x_0$. 
In addition to point-sources, \cite{fomel2009} also used factoring to treat ``plane-wave sources'' (i.e., constant Dirichlet boundary conditions specified on a straight line in 2D). In the latter case, there is no singularities in the solutions at the boundary (so, the numerics for the original/unfactored version is still first-order accurate), but a factored version still yields lower errors in many examples.

The original factoring in \cite{fomel2009} used an ansatz $u(\x) = T(\x)\tau(\x)$ to derive a new factored PDE for $\tau$,
which was then solved with the boundary condition $\tau(\x_0)=1.$
In this paper we rely on a slightly simpler additive splitting\footnote{
Throughout the paper we refer to this approach as ``additive factoring'' to stay consistent
with the terminology used in prior literature.} 
introduced in \cite{luo2012fast} : assuming that
$u(\x) = T(\x)+ \tau(\x),$ we find $\tau$ by solving
\begin{equation} 
\label{eq:factor}
\left| \nabla T(\x) + \nabla\tau(\x) \right| F(\x) \, = \, 1, 
\end{equation}
with the boundary condition $\tau(\x_0) = 0.$
Similarly to \eqref{eq:Eikonal} this can be discretized using upwind finite differences and then solved on $\Omega$ (see section \ref{ss:modified}).
In all our convergence figures we refer to this approach as the ``global factoring'' (shown by a solid blue line). 
In Figures \ref{f:nonobst_1}-\ref{f:nonobst_3} it is clear that global factoring has a clean linear convergence. 
Moreover, in Figure \ref{f:nonobst_1} it exhibits smaller errors than either the original/``unfactored'' discretization or even the Lagrangian-initialized version regardless of the grid resolution.  However, as examples in Figures \ref{f:nonobst_2} and \ref{f:nonobst_3} show, on relatively coarse grids it can be actually less accurate than the unfactored discretization -- particularly when the characteristics are far from straight lines.  This phenomenon has not been examined in prior literature, 
but it is hardly surprising: farther from the point source, $u$ and $T$ can be quite different, and if the second derivatives of $u$ are smaller, this will result in larger local truncation errors when computing $\tau$.  

We further examine a ``localized factoring'' version of this idea.  For small $r,$ this could be posed as a 2-stage process: first solve \eqref{eq:factor} on $B$ and then switch to solving \eqref{eq:Eikonal} on $\Omega \backslash B.$  However, we have found that another interpretation is more suitable, particularly when characteristics are highly curved: solve  \eqref{eq:factor} on the entire $\Omega$ but defining 
$T(\x)=0$ on $\Omega \backslash B.$  We note that several recent papers have already considered such hybrid/localized factoring motivated  by decreasing the computational cost (since $\nabla T = 0$ on most of the domain) \cite{noble2014} and by the need for additional properties of $T$ when pursuing  a higher-order accurate discretization \cite{luo2014high}. 
Here, however, we show that the localized factoring (shown by a blue dashed lines on convergence plots) has some accuracy advantages even with the first-order upwind discretization.  In our first example (Figure \ref{f:nonobst_1}), $u \approx T$ remains true on the whole $\Omega$ and characteristics are fairly close to straight lines; so, the global factoring is more accurate.
But in Figure \ref{f:nonobst_2} this is no longer the case, and the localized factoring is clearly preferable.

Localized factoring is also often advantageous (and more natural) when dealing with multiple point sources.
Consider, for example, the speed function specified in \eqref{eq:single_src} with $s_0 = 0.5, \, \bv = (5,20)$ and two point sources:
$\x_0 = (0,0)$ and $\x_1=(0,0.8).$   (See Figure \ref{f:nonobst_3}.)  If $Q = \{\x_0\}$ or $Q = \{\x_1\}$, the respective value functions $u_0$ and $u_1$ 
are specified by formula \eqref{eq:single_src}.  For the two point-sources case, $Q = \{\x_0, \, \x_1\}$, the value function $u(\x) = \min \left(u_0(\x), \, u_1(\x)\right)$ is no longer smooth: $\nabla u$ is undefined at the points from which the optimal paths to $\x_0$  and $\x_1$ are equally good.
Since the characteristics run into the shockline rather than originate from it, this does not degrade the rate of convergence, but 
the rarefaction fans remain a challenge.
In global factoring, there is a number of choices to capture in $T$ the singularities at both point sources.  We use
\[ T_0(\x) = \frac{\left| \x - \x_0 \right|}{F(\x_0)},\qquad T_1(\x) = \frac{\left| \x - \x_1 \right|}{F(\x_1)},\qquad T(\x) = \min\big(T_0(\x),T_1(\x)\big), \] 
but the convergence results based on $T= T_0 + T_1$ are qualitatively similar (though the $L^1$ error becomes significantly larger).  
For the localized factoring, we simply set $T(\x)=T_i(\x)$ on $B_r(x_i)$
and $T(\x)=0$ everywhere else.
Figure \ref{f:nonobst_3} shows that both versions of factoring exhibit linear convergence, but the localized factoring yields significantly smaller errors both in  $L^{\infty}$ and $L^1$ norms.

\subsection{Fast methods for an (unfactored) Eikonal}
\label{ss:fastEik}

Since an Eikonal equation arises in so many applications, there has been a number of fast numerical methods developed for it in the last 20 years.  Most of these methods mirror the logic of classical label-setting and label-correcting algorithms for finding shortest paths on graphs \cite{bertsekas1998networks}.  In this discrete setting, an equation for the min-time-to-exit $U_i$ starting from each node $\x_i$ is 
posed using the min-time-to-exit values ($U_j$'s) at its neighboring nodes ($\x_j$'s).
Dijkstra's method \cite{Dijkstra} is perhaps the most famous of label-setting algorithms for graphs with positive edge-weights.
It is based on the idea of {\em monotone causality}: an optimal path from $x_i$ starts with a transition to some neighboring node $\x_{j^*}$ and, since 
all edge weights are positive, this implies $U_i > U_{j^*}.$  Thus,  even if we don't know $\x_{j^*}$ a priori, $U_i$ can be still computed based 
on the set of all {\em smaller} neighboring values.  
Dijkstra's method exploits this observation to dynamically uncover the correct node-ordering, decoupling the system of equations for $U_i$'s.
The nodes are split into three sets: Accepted (with permanently fixed $U$ values), Considered (with tentatively computed $U$ values) and Far (with $U$ values assumed to be $+\infty$).
At each stage of the algorithm, the smallest of Considered $U$ values is declared Accepted and the values at its immediate not-yet-Accepted neighbors are recomputed.  On a graph with $M$ nodes and bounded node connectivity, this yields the overall complexity of $O(M \log M)$ due to the use of a heap to sort the Considered values.  
An additional useful feature of this approach is that the $U$ values are fixed/accepted after a small number of updates on an (incrementally growing) part of the graph.  
Many label-correcting algorithms aim to mimic this property, but without using expensive data structures to sort any of the values.  
Unlike in label-setting algorithms, they cannot provide an a priori upper bound on the number of times each $U_i$ might be updated.  As a result, their worst-case complexity is $O(M^2),$ but on many types of graphs their average-case behavior has been observed to be at least as good as that of label-setting techniques \cite{bertsekas1998networks}.

For Eikonal PDEs discretized on Cartesian grids, the Dijkstra-like approach was first introduced in 
Tsitsiklis's Algorithm \cite{tstitsiklis1996} and Sethian's Fast Marching Method (FMM) \cite{sethian1996fast}.
The latter was further extended to simplicial mesh discretizations  in $R^n$ and on manifolds \cite{KimmSethTria,SethVlad1}, 
to higher-order accurate numerical schemes \cite{SethSIAM, SethVlad1}, and to Hamilton-Jacobi-Bellman equations, which can be viewed as anisotropic generalizations of the Eikonal \cite{SethVlad2, SethVlad3, AltonMitchell2, Mirebeau3, DahiyaCameron2017}.
The major difficulty in applying Dijkstra-like ideas in continuous setting is that, unlike in problems on graphs,  a value of $u$ at a gridpoint $\x_{i,j}$
will depend on several neighboring values used to approximate $\nabla u(\x_{i,j}).$  To obtain the same monotone causality, this $u(\x_{i,j})$ has to be 
always larger than {\em all} of such contributing $u$ values adjacent to $\x_{i,j}.$  This property is enjoyed by only some discretiations of \eqref{eq:Eikonal}, including the first-order accurate upwind scheme \eqref{eq:scheme}: a direct verification shows that $u_x$ and $u_y$ are never approximated using any neighbors larger that $u_{i,j}.$  

Another popular class of efficient Eikonal solvers is Fast Sweeping Methods (FSM) \cite{zhao2005fast}, which solve the system \eqref{eq:scheme} by Gauss-Seidel iterations, but changing the order in which the gridpoints are updated from iteration to iteration.  When the direction of the current sweep is aligned with the general direction of characteristics, many gridpoints will receive correct values in a single ``sweep''.  If marching-type methods attempt to uncover the correct gridpoint-ordering dynamically, in FSM the idea is to alternate through a number of geometrically motivated orderings. (In 2D problems: from northeast, from northwest, from southwest, and from southeast.)  The resulting Eikonal solvers have $O(M)$ complexity, but with a hidden constant factor, which depends on $F$ and the grid orientation, and cannot be bound a priori. 
These techniques have also been extended to anisotropic (and even non-convex) Hamilton-Jacobi equations \cite{TsaiChengOsherZhao, Kao2004}, 
as well as higher-order finite-difference (e.g., \cite{Zhang2006}) and discontinuous Galerkin (e.g., \cite{Zhang2011}) discretizations.
Hybrid two-scale methods, combining the best features of marching and sweeping, were more recently introduced in \cite{ChacVlad1, ChacVlad2}.
We also refer to \cite{ChacVlad1} for a comprehensive review of other fast solvers inspired by label-correcting algorithms.

\subsection{Modified Fast Marching for factored Eikonal}
\label{ss:modified}
We start by simplifying the original upwind discretization scheme \eqref{eq:scheme} for the unfactored Eikonal.
Focusing on a gridpoint $\x_{i,j}$ we define its smallest horizontal and vertical neighboring values:
$\uH = \min( u_{i-1,j}, \, u_{i+1,j})$ and $\uV = \min( u_{i,j-1}, \, u_{i,j+1}).$  
If both of these values are needed to compute $u_{i,j},$
then \eqref{eq:scheme} becomes equivalent to a quadratic equation
$(u_{i,j} - \uH)^2 \, + \, (u_{i,j} - \uV)^2 \, = \, h^2 / F^2_{i,j}.$
We are interested in its smallest real root satisfying an {\em upwinding condition} $u_{i,j} \geq \max(\uH, \, \uV).$  
If there is no root satisfying it, then $u_{i,j}$ should instead be computed from a {\em one-sided-update}
$u_{i,j} \, = \, \min(\uH, \, \uV) + h/F_{i,j},$ which corresponds to the case where 
either $\max\{D_{i,j}^{-x}u, -D_{i,j}^{+x}u, 0\}$ or $\max\{D_{i,j}^{-y}u, -D_{i,j}^{+y}u, 0\}$
evaluates to zero.  This procedure is monotone causal by construction and its equivalence to \eqref{eq:scheme}
was demonstrated in \cite{sethian1996fast}, making a Dijkstra-like computational approach suitable.

Recalling the ``additive factoring'' ansatz $u(\x) = T(\x) + \tau(\x),$ we now define the upwind vertical and horizontal neighboring values
for the new unknown $\tau_{i,j}$ but basing the comparison on $u$ rather than on $\tau$ itself
and using the flags $\kH, \kV \in \{-1,\, 1\}$ to identify the selected neighbors.  
More specifically,
\begin{equation*}
\begin{cases}
\tauH = \tau_{i-1,j} \text{ and } \kH = 1, &
\text{if } (T_{i-1,j} + \tau_{i-1,j}) \, <  \, (T_{i+1,j} + \tau_{i+1,j});\\ 
\tauH = \tau_{i+1,j} \text{ and } \kH = -1, &
\text{otherwise;}
\end{cases}
\end{equation*}
with $(\tauV, \, \kV)$ similarly defined based on the vertical neighbors.
Since the partial derivatives of $T$ are known, the corresponding quadratic equation is
\begin{equation}
\label{eq:quadratic}
\left( \kH \frac{\partial T_{i,j}}{\partial x} + \frac{\tau_{i,j}-\tauH}{h} \right)^2 
\, + \,
\left( \kV \frac{\partial T_{i,j}}{\partial y} + \frac{\tau_{i,j}-\tau_{v}}{h} \right)^2 
\; = \; \frac{1}{F_{i,j}^2}.
\end{equation}
We are interested in its smallest real root satisfying a similarly modified {\em upwinding condition} 
\begin{equation}
\label{eq:causality}
T_{i,j} + \tau_{i,j} \; \geq \; 
%\max(\uH, \, \uV ) \, = \, 
\max\left(
\min_{k=\pm1} \left\{T_{i+k,j} + \tau_{i+k,j}\right\},
\,
\min_{k=\pm1} \left\{T_{i,j+k} + \tau_{i,j+k}\right\}
\right).
\end{equation}
If there is no such root, then $\tau_{i,j}$ should instead be computed from a {\em one-sided-update} as the smaller of the two values 
corresponding to
$$
\kH \frac{\partial T_{i,j}}{\partial x} + \frac{\tau_{i,j}-\tauH}{h} \, = \, \frac{1}{F_{i,j}}, 
\qquad 
\text{ and }
\qquad
\kV \frac{\partial T_{i,j}}{\partial y} + \frac{\tau_{i,j}-\tauV}{h} \, = \, \frac{1}{F_{i,j}}.
$$
In other words,
\begin{equation}
\label{eq:oneside2}
\tau_{i,j} \; = \;  \min \left\{ 
\tauH - h \kH \frac{\partial T_{i,j}}{\partial x}, \;
\tauV - h \kV \frac{\partial T_{i,j}}{\partial y}
\right\}
\, + \, \frac{h}{F_{i,j}}.
\end{equation}
The above is a full recipe for computing $\tau_{i,j}$ if all of the neighboring grid values are already known.
But since $\tau$ is a priori known only on $Q,$ this yields a large coupled system of discretized equations
(one per each gridpoint in $\Omega \backslash Q$).  This system was first treated iteratively via Fast Sweeping \cite{luo2012fast},
but the monotone causality makes a Dijkstra-like approach applicable as well.
As in the original Fast Marching, the Considered gridpoints are sorted using a binary heap with a $O(M \log M)$ computational complexity;
however, the sorting criterion is based on $u$ rather than on $\tau$ values.  
The resulting method is summarized in Algorithm \ref{alg:mfmm}.  It is very similar to a modified FMM recently introduced for the case of ``multiplicative factoring'' in \cite{treister2016fast}.   
%It becomes equivalent to the original FMM if we use $T=0.$

\begin{algorithm}[!htb]
\caption{Modified Fast Marching Method}
\label{alg:mfmm}
\KwIn{source point $\x_0$, speed function $F(\x)$}
%\KwOut{}
$\empty$\\
Initialize $\tau(\x_0) := 0$ and $\tau(\x) := +\infty$ for all gridpoints $\x \neq \x_0$\;
Initialize Considered $:= \{\x_0\}$ and Accepted $:= \emptyset$\;
$\empty$\\
\While{Considered $\neq \emptyset$}{
	Find a Considered gridpoint  $\xhat$ whose $(T+\tau)$ value is the smallest\;
	Move $\xhat$ to Accepted list\;
$\empty$\\	
	\For{all not-yet-Accepted neighbors $\x_{i,j}$ of $\xhat$}{
        Find $\tauH, \tauV$ and solve the quadratic equation \eqref{eq:quadratic} for $\tau_{i,j}^{new}$\;
		\If{$\tau_{i,j}^{new}$ does not satisfy the upwinding condition \eqref{eq:causality}}{
			Use a one-sided update formula \eqref{eq:oneside2} to (re-)compute $\tau_{i,j}^{new}$\;
		}
		\If{$\tau_{i,j}^{new} \, < \, \tau_{i,j}$}{
			Set $\tau_{i,j} \, := \, \tau_{i,j}^{new}$\;
		}
	}
}
\end{algorithm}

%-----------------------------------------------------------------------

% !TEX root = local_factoring_for_Journal.tex

\section{Rarefaction fans at obstacle corners}
\label{s:bad_corners}

\begin{figure}[!htb]
\centering
\subfigure[multiple obstacles]{
\label{f:tra1}
\includegraphics[width=0.48\textwidth]{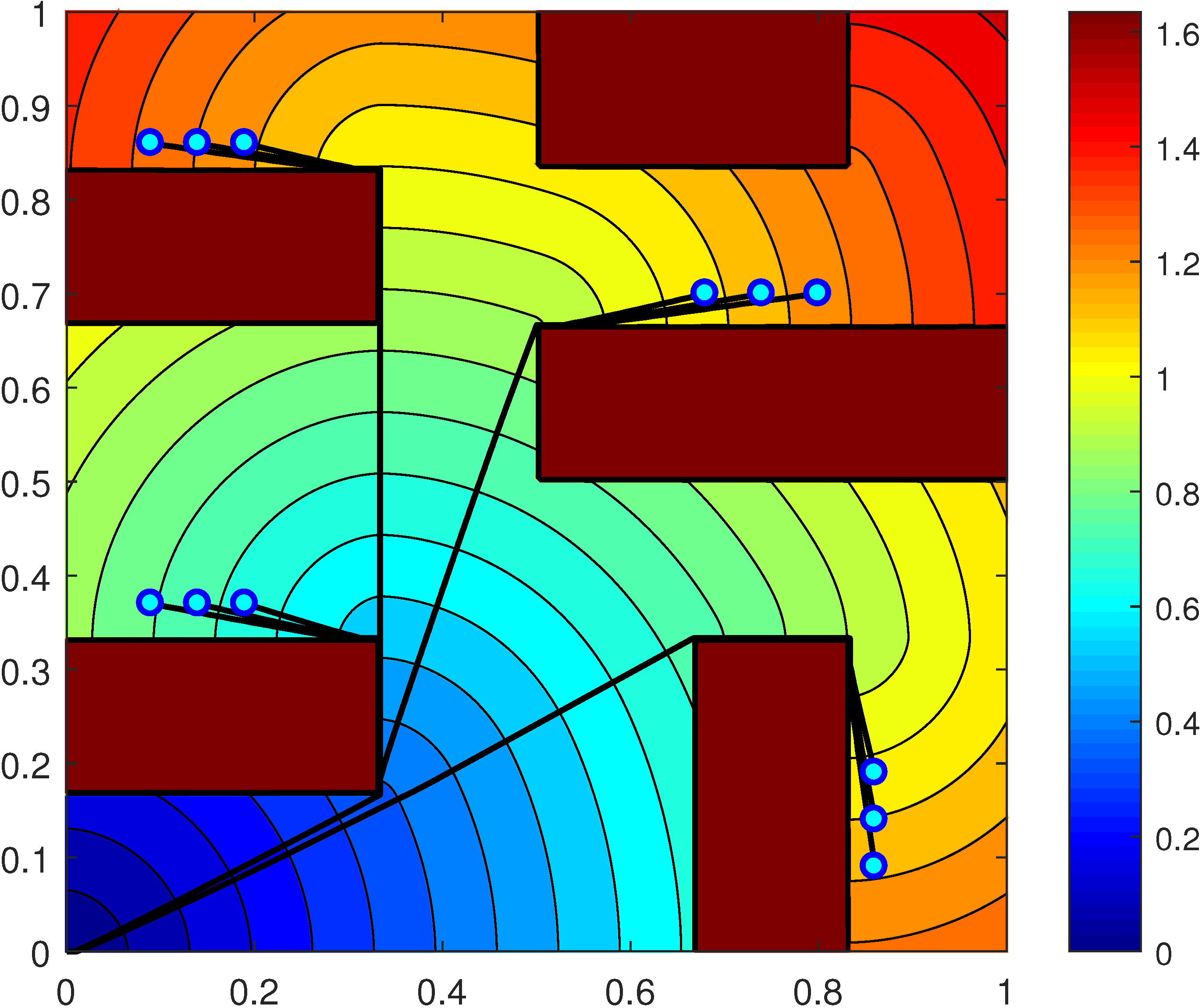}}
\subfigure[distorted trajectories]{
\label{f:tra2}
\includegraphics[width=0.40\textwidth]{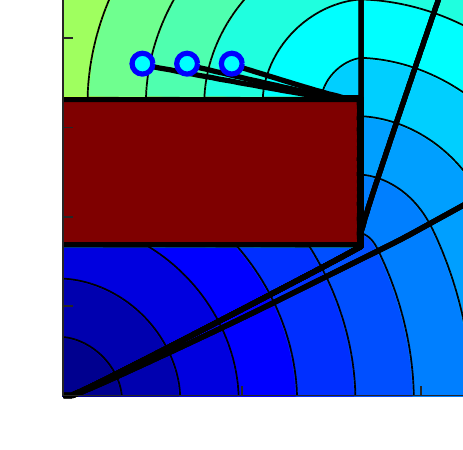}}
\caption{A maze navigation example: non-permiable obstacles with $F=1$ on the rest of $\Omega.$
(Left) The level sets of the value function $u$ computed by the Fast Marching Method on a $240 \times 240$ grid
and approximate optimal trajectories to the origin from 12 starting locations.  (Right) A zoomed version to highlight the incorrect direction of ``optimal'' trajectories in the rarefaction fans at obstacle corners. }
\label{f:tra}
\end{figure}

Even though all prior work on factored Eikonal equation was focused on isolated point sources, there are other well-known situations where rarefaction fans can arise.  As a simple example in Figure \ref{f:hessian} shows, they can easily develop at the corners of obstacles (which are viewed as a part 
of $R^2 \backslash \Omega$)
 or, more generally at any points on $\partial \Omega$ where the boundary is non-smooth and the interior angle is larger than $\pi.$
 These non-point-source rarefaction fans result in a similar degradation of convergence rate for standard numerical methods and also lead to unpleasant artifacts in optimal trajectory approximations obtained by following $(-\nabla u$) to the target set $Q$.  Figure \ref{f:tra} shows several such trajectories in a ``maze navigation'' problem.  All of these trajectories should be piecewise-linear, with their directions only changing at obstacle corners.  A zoomed version in Figure \ref{f:tra2} clearly shows that they often approach an obstacle too early, following its boundary to the corner and yielding longer paths.  
Similar artifacts are common in determining parts of the domain visible by an observer \cite{SethAdal1997} and in multiobjective path-planning \cite{KumarVlad, MitchellSastry_Multio_Published}. 
A natural question (and the focus of this paper) is whether factoring techniques
can be used to alleviate this problem.
In section \ref{ss:nonpermeable_examples} we demonstrate experimentally that the ``global factoring'' is not suitable for this task.  
On the other hand, the localized factoring works, but adopting it to corner-induced rarefaction fans presents two new challenges.
First of all,  not all obstacle corners produce this effect; e.g., see the lower left corner in Figure  \ref{f:hessian_v}.
\begin{defn}
\label{def:badcorner}
An obstacle corner $\xtilde$ is ``\textit{regular}'' if the characteristic leading to it from $Q$ points {\em into} that obstacle.
(I.e., if an optimal trajectory starts from $\xtilde$ in the direction $\ba$, 
then $(-\ba)$ should point into the obstacle.)
An obstacle corner is ``\textit{rarefying}'' if it is not regular.
\end{defn}
So, even though the rarefying corners are not known in advance, we can identify them dynamically, checking the above condition when the corresponding corner $\xtilde$ becomes Accepted in Fast Marching Method 
and approximating the optimal 
$\ba = \frac{-\nabla u(\xtilde)}{\left| \nabla u(\xtilde) \right|}$ using $\xtilde$'s previously Accepted upwind neighbors.  
The resulting {\em ``just-in-time localized factoring''} method is detailed in Algorithm \ref{alg:dynamic}.
It maintains a list of identified rarefaction fans, with each entry containing the center of the fan (either a point source or a rarefying corner) and the corresponding localized function $T$.  
The algorithm is formulated in terms of $u$, with the corresponding $\tau$ values computed on the fly once the appropriate localized $T$ is selected.

\begin{algorithm}[!htb]
\caption{Just-in-time Localized Factoring}
\label{alg:dynamic}
\KwIn{source point $\x_0$, speed function $F(\x)$, fixed radius $r$}
%\KwOut{}
$\empty$\\
Initialize $u(\x_0) := 0$ and $u(\x) := +\infty$ for all gridpoints $\x \neq \x_0$\;
Initialize Considered $:= \{\x_0\}$ and Accepted $:= \emptyset$\;
Initialize FanList $:= \left\{ \left( \x_0, \, T^{\x_0} = |\x-\x_0|/F(\x_0)  \right) \right\}$\;
$\empty$\\
\While{Considered $\neq \emptyset$}{
	Find a Considered gridpoint  $\xhat$ whose $u$ value is the smallest\;
	Move $\xhat$ to Accepted list\;
	\If{$\xhat$ is a rarefying corner}{
		Build a suitable $T^{\xhat}$ using formula \eqref{eq:factorT}\;
		Add an entry $\left( \xhat, T^{\xhat} \right)$ to FanList\;
	}
	
$\empty$\\
         \For{all not-yet-Accepted neighbors $\x_{i,j}$ of $\xhat$}{
         	Check if $\x_{i,j}$ is within distance $r$ from any $\xtilde$ on FanList\\
		and identify the appropriate $T=T^{\xtilde}$ (use $T=0$ by default)\;
		Given the current $u$ values at $\x_{i,j}$ \& its neighbors,\\
		define their $\tau$ values as $\tau = (u - T)$ based on $T$ selected for $\x_{i,j}$\;  
$\empty$\\
	        Find $\tauH, \tauV$ and solve the quadratic equation \eqref{eq:quadratic} for $\tau_{i,j}^{new}$\;
		\If{$\tau_{i,j}^{new}$ does not satisfy the upwinding condition \eqref{eq:causality}}{
			Use a one-sided update formula \eqref{eq:oneside2} to (re-)compute $\tau_{i,j}^{new}$\;
		}
		$u_{i,j}^{new} \, := \, \tau_{i,j}^{new} + T_{i,j}$\;
		\If{$u^{new} \, < \, u_{i,j}$}{
			Set $u_{i,j} \, := \, u_{i,j}^{new}$\;
		}
	}
}
\end{algorithm}

\begin{figure}[!htb]
\centering
\subfigure[level sets of $u$]{
\label{f:illustration_1}
\iftoggle{ForJournal}{\includegraphics[scale=.34,clip]{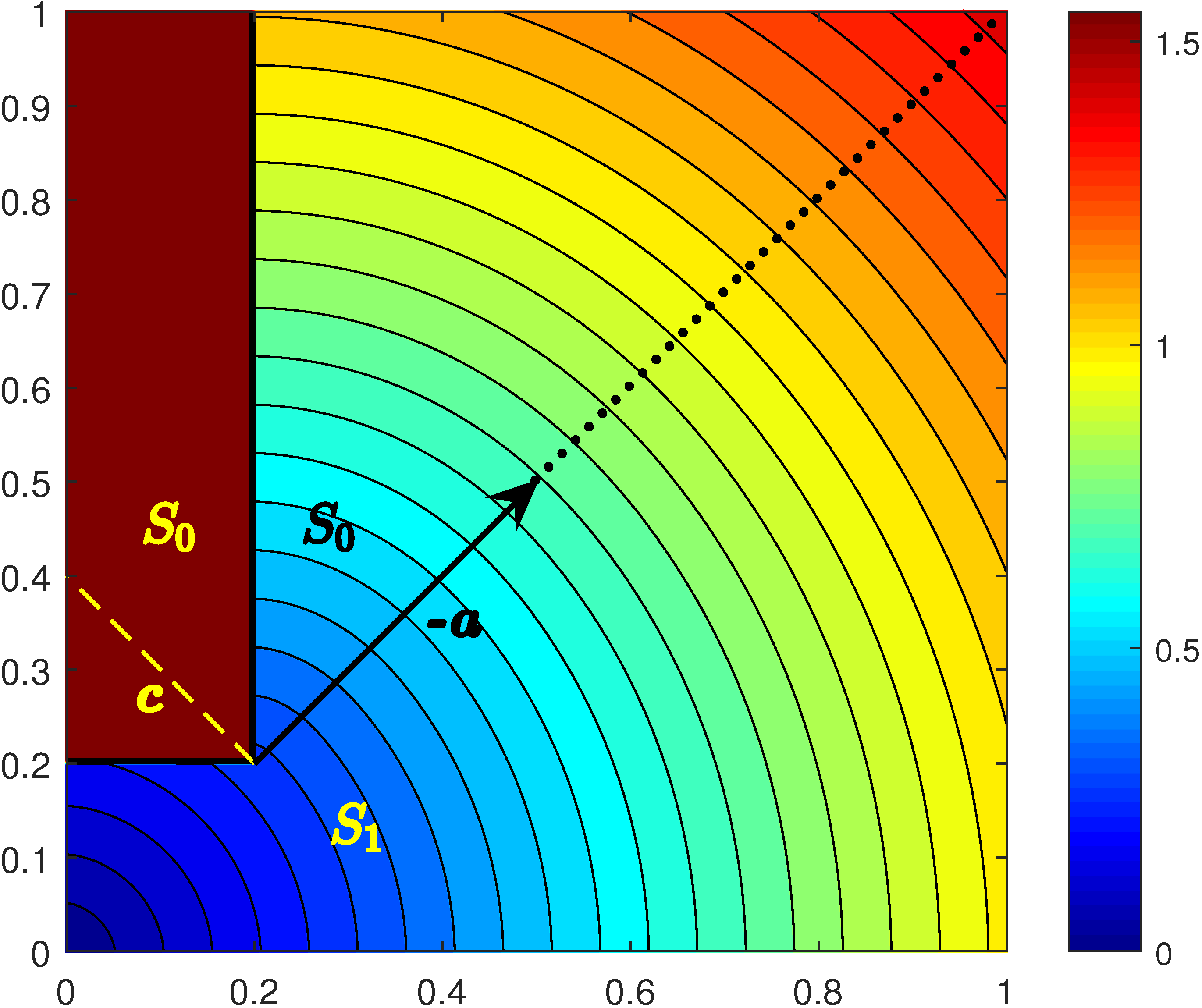}}{\includegraphics[scale=.43,clip]{illustration_1}}}
\subfigure[domain splitting]{
\label{f:illustration_2}
\iftoggle{ForJournal}{\includegraphics[scale=.34,clip]{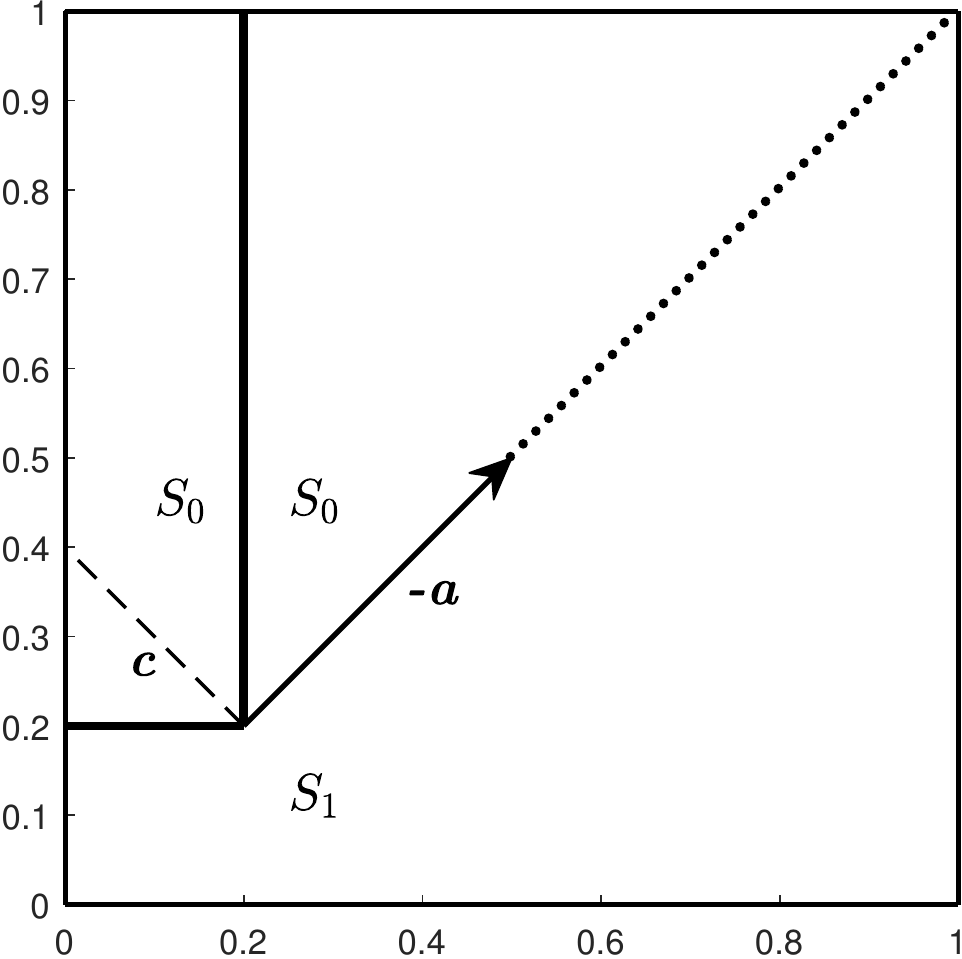}}{\includegraphics[scale=.43,clip]{regions}}}
\subfigure[level sets of $T$]{
\label{f:illustration_3}
\iftoggle{ForJournal}{\includegraphics[scale=.34,clip]{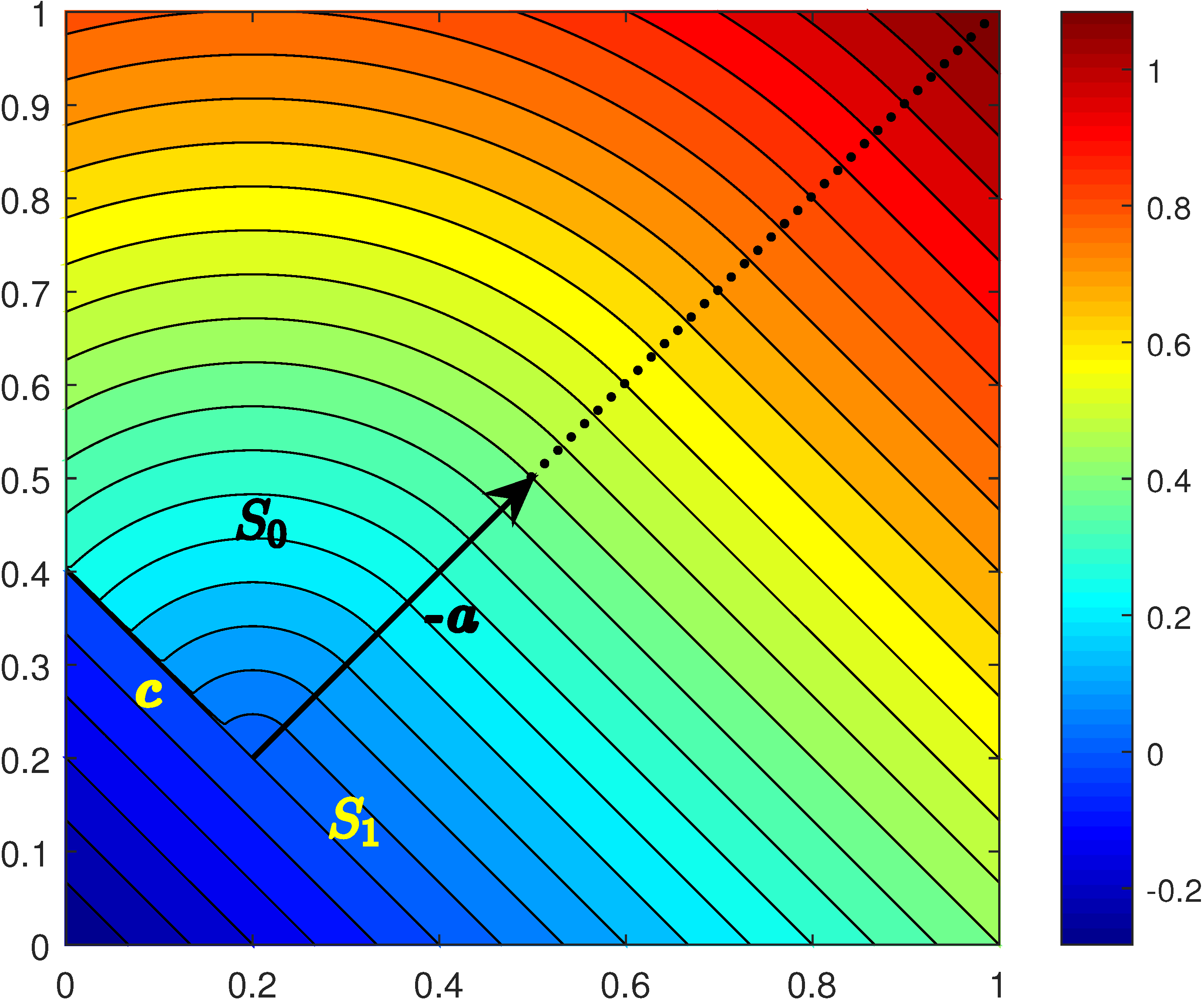}}{\includegraphics[scale=.43,clip]{t_value}}}
\caption{A simple example with one rarefying corner. 
(Left) level curves of the domain-restricted distance to a point source.
(Center) a dynamic domain splitting based on the rarefaction fan.  
(Right) the level sets of a ``cone+plane'' function $T$ capturing the correct rarefaction behavior.}
\label{f:illustration}
\end{figure}

The second difficulty is to define a suitable $T$ that will be used for factoring when updating all not-yet-Accepted gridpoints in  $B_r(\xtilde).$ 
Intuitively, it might seem that a cone-like $T = \frac{\left| \x - \xtilde \right|}{F(\xtilde)}$ is the right choice, similarly to our handling of point sources.
However, as we show in section \ref{ss:nonpermeable_examples}, this choice does not yield the desired rate of convergence.  
This is due to the fact that such corner-born rarefaction fans are not radially symmetric.
They only exist for $u > u(\xtilde)$ in the sector between a part of obstacle boundary and the characteristic passing through $\xtilde$;
see Figure \ref{f:hessian} and an even simpler example in Figure \ref{f:illustration}. 
Note that, outside of that rarefaction sector, the second derivatives of $u$ are bounded.
But using the ansatz $u = T + \tau$ with the above cone-like $T$ would introduce unbounded second derivatives in $\tau$ 
%on $S_1 \cap B_r(\xtilde)$, 
on non-rarefying parts of $B_r(\xtilde)$, thus degrading the rate of convergence. 
%on the entire domain.  
Therefore, we need to construct $T$ which is cone-like only in the correct sector and remains smooth on the entire not-yet-Accepted portion of the domain.  This yields a ``cone+plane'' version of $T$ shown in Figure \ref{f:illustration_3}. 
Assuming that $\bc$ is a unit vector bisecting the obstacle corner at $\xtilde,$ we can now split the plane into two sets
$$S_0 = \left\{ \x \in \Omega \mid (\x -\xtilde) \text{ is between $\bc$ and $(-\ba)$ }\right\};
\qquad \qquad
S_1 = \Omega \backslash S_0.$$
Analytically, we can define $T$ as follows
\begin{equation}
\label{eq:factorT}
T(\x) \; =\;
\begin{dcases}
\dfrac{\left| \x-\xtilde \right|}{F(\xtilde)}, \quad & \x \in S_{0} \\
%\\
\dfrac{-\ba \cdot (\x-\xtilde)}{F(\xtilde)}, \quad & \x \in S_{1}.
\end{dcases}
\end{equation}
The resulting $T$ is not continuous along $\bc$, but this will not matter since the discontinuity is hidden within the obstacle. 
The gradient of $T$ is also continuous wherever $T$ is, though the second derivatives are bounded but discontinuous along $(-\ba)$.
Our numerical results  show that this $T$ fully recovers the first-order convergence of the numerical solution.

\begin{rem} \label{rem:discont}
The bisector of obstacle corner is not the only choice for $\bc$.  Any directions falling inside the obstacle will work just as well since the idea is to ``hide'' the discontinuity line of $T.$ 
For rectangular obstacles, it might feel more natural to choose $\bc$ orthogonal to the characteristic direction $\ba.$
However, we prefer the bisector simply because it is a safe choice for arbitrary polygonal obstacles,
which can be handled by a version of Algorithm \ref{alg:dynamic} on (obstacle-fitted) triangulated meshes.

Another possibility is to hide the $T$'s discontinuity in the part of $\Omega$ already Accepted (e.g., along $\ba$) by the time this rarefying corner is identified.  This is the approach we use in section \ref{ss:discont}, when dealing with ``slowly permeable obstacles''.%, inside  which the value function is defined and finite.
\end{rem}

\begin{rem} \label{rem:in_S0_only}
Our approach uses local factoring with a continuous $T$ 
%defined by formula \eqref{eq:factorT} 
on the entire $B_r(\xtilde) \cap \Omega$.  A variant of the same idea is to  employ local factoring with a standard cone-like $T = \frac{\left| \x - \xtilde \right|}{F(\xtilde)}$ but only when updating gridpoints from a subset  $B_r(\xtilde)  \cap \Omega \cap S_0$.  (The difference is that one would use $T=0$ when updating gridpoints in $B_r(\xtilde)  \cap \Omega \cap S_1$.)  While we do not formally include this variant in our convergence studies in section \ref{ss:nonpermeable_examples}, its performance appears to be quite similar.  E.g.,  for the above example from Figure \ref{f:illustration}, the alternative version also shows the first-order convergence, but with $L^{\infty}$ errors $\approx10\%$ larger than those resulting from the ``cone+plane'' formula \eqref{eq:factorT}.
\end{rem}

\vspace*{2mm}

\noindent
We close this section by discussing a subtle property implicitly used in our approach.
The above construction relies on having a sufficiently accurate representation of the characteristic direction $\ba$ at each rarefying corner.  This might seem unreasonable: if our numerical approximation of $u$ is only $O(h)$ accurate (as is the case in formulas \eqref{eq:scheme} and (\ref{eq:quadratic}-\ref{eq:oneside2})), then one could expect the resulting finite difference approximation of $\nabla u$ to be completely inaccurate.   The same argument would imply that optimal trajectories also cannot be reliably approximated based on any first-order accurate representation of the value function.  However, there is ample experimental evidence that such trajectory approximation works in practice.  See, for example, the optimal trajectories in Figure \ref{f:tra},\ref{f:maze_1},\ref{f:maze_2},\ref{f:perm_1} and in many optimal control and seismic imaging publications with and without factoring.  
The fact that this gradient approximation is in fact $O(h)$ accurate is also confirmed by a numerical study in \cite{BenamouLuoZhao} and is instrumental for constructing other (compact stencil, second-order) schemes for the Eikonal \cite{SethVlad1, BenamouLuoZhao}.  

A plausible explanation for this ``superconvergence'' phenomenon is that the error in $u$-approximation is sufficiently ``smooth'', resulting in a convergent $\nabla u$--approximation despite the use of divided differences.   To the best of our knowledge, this property has not been proven for general Eikonal PDEs, though it has been rigorously demonstrated for the distance-to-a-point computations and for constant coefficient linear advection equations\cite[Appendix B]{SLuoThesis}.
In our current context, we use the same idea to conjecture that the $\ba$-dependent approximation of the localized $T$ is sufficiently accurate to recover the full first-order accuracy in $u$ computations with dynamic factoring.  This conjecture appears to be fully confirmed by the convergence rates observed in our numerical experiments throughout this paper.

%-----------------------------------------------------------------------

% !TEX root = local_factoring_for_Journal.tex

\subsection{Numerical Examples} 
\label{ss:nonpermeable_examples}
As a first numerical test, we consider a simple example from Figure \ref{f:illustration_1}: a domain-constrained distance $u$ to the origin
on $\Omega=[0,1]\times[0,1] \backslash \Omega_{ob}$, with a single rectangular obstacle 
$\Omega_{ob}=(0,0.2)\times(0.2,1).$ 
Since $F=1,$ all optimal paths are piecewise-linear.  
According to Definition \ref{def:badcorner}, we find that the corner at $\xtilde = (0.2,0.2)$ is rarefying.
As a result, the problem contains two rarefaction fans: one at this corner and the other at the point source $\x_0 = (0,0).$

\begin{figure}[!htb]
\centering
\subfigure[level sets]{
\label{f:obst_1_v}
\includegraphics[width=0.44\textwidth,clip]{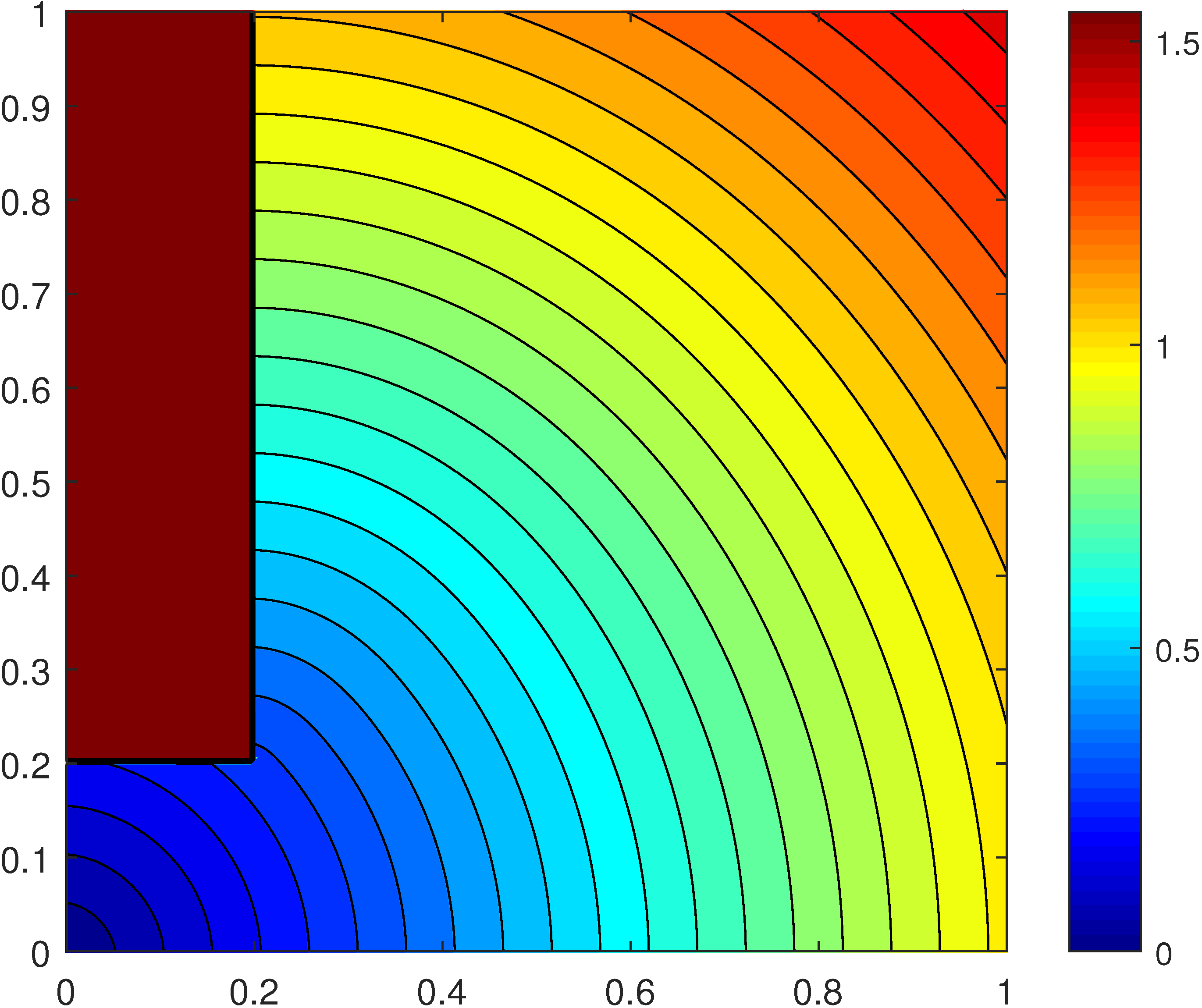}}
\subfigure[$L^{\infty}$ error]{
\label{f:obst_1_c}
\iftoggle{ForJournal}{\includegraphics[scale=.42,clip]{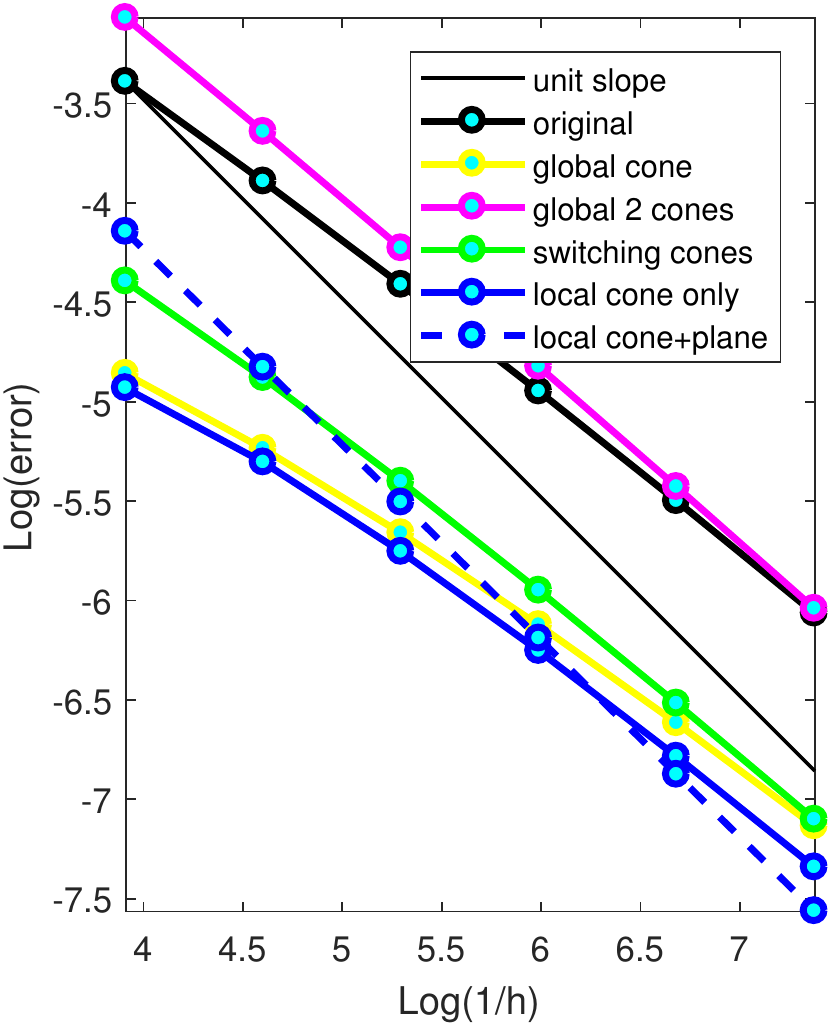}}{\includegraphics[scale=.53,clip]{obst_1_c}}}
\caption{Convergence of several methods for a simple obstacle case of Figure \ref{f:illustration_1}. 
}
\label{f:obst_1}
\end{figure}

\noindent
We test the accuracy of several methods:
\begin{enumerate}
\item \textbf{original:} Original (un-factored) Eikonal solved on the entire $\Omega$ with the original Fast Marching Method.

\item \textbf{global cone:} Global factoring using $T(\x) = \left| \x-\x_0 \right| / F(\x_0)$ on the entire $\Omega$ with Algorithm \ref{alg:mfmm}.

\item \textbf{global 2 cones:} Global factoring using $T(\x) = \left| \x-\x_0 \right| / F(\x_0) \, + \, \left| \x-\xtilde \right| / F(\xtilde)$ on the entire $\Omega$ with Algorithm \ref{alg:mfmm}.

\item \textbf{switching cones:} Global factoring using $T(\x) = \left| \x-\x_0 \right| / F(\x_0)$ until $\xtilde$ is accepted and then switch to
global factoring using $T(\x) = \left| \x-\xtilde \right| / F(\xtilde)$ on the rest of $\Omega.$

\item \textbf{localized cone only:} Just-in-time localized factoring Algorithm \ref{alg:dynamic} with $T(\x) = \left| \x-\x_0 \right| / F(\x_0)$ on $B_r(\x_0)$ and 
another cone-like $T(\x) = \left| \x-\xtilde \right| / F(\xtilde)$ on $B_r(\xtilde)$.

\item \textbf{localized cone+plane:} Just-in-time localized factoring Algorithm \ref{alg:dynamic} with $T(\x) = \left| \x-\x_0 \right| / F(\x_0)$ on $B_r(\x_0)$ and 
a dynamically defined ``cone+plane'' $T(\x)$ specified by formula \eqref{eq:factorT} on $B_r(\xtilde)$.
\end{enumerate}
The first four of these are included to show that the corner-induced rarefaction fans do indeed degrade the rate of convergence
and the issue cannot be addressed by global factoring.  
Accuracy of all methods is tested using a range of gridsizes ($h = \tfrac{1}{50}2^{-k}$, where $k=0, \ldots, 5$).
The localized factoring is based on $r=0.18$.
As Figure \ref{f:obst_1_c} clearly shows, only the last method actually exhibits the first-order of convergence. 
Even though the usual global factoring (method 2) starts out with smaller errors on coarser meshes, 
it becomes worse than our preferred approach (method 6) for smaller values of $h$.    
The fact that method 5 has a similar performance degradation 
proves the importance  of choosing the correct localized factoring function $T$.
 
\begin{figure}[!htb]
\centering
\subfigure[maze: original]{
\label{f:maze_1_v1}
\includegraphics[width=0.47\textwidth,clip]{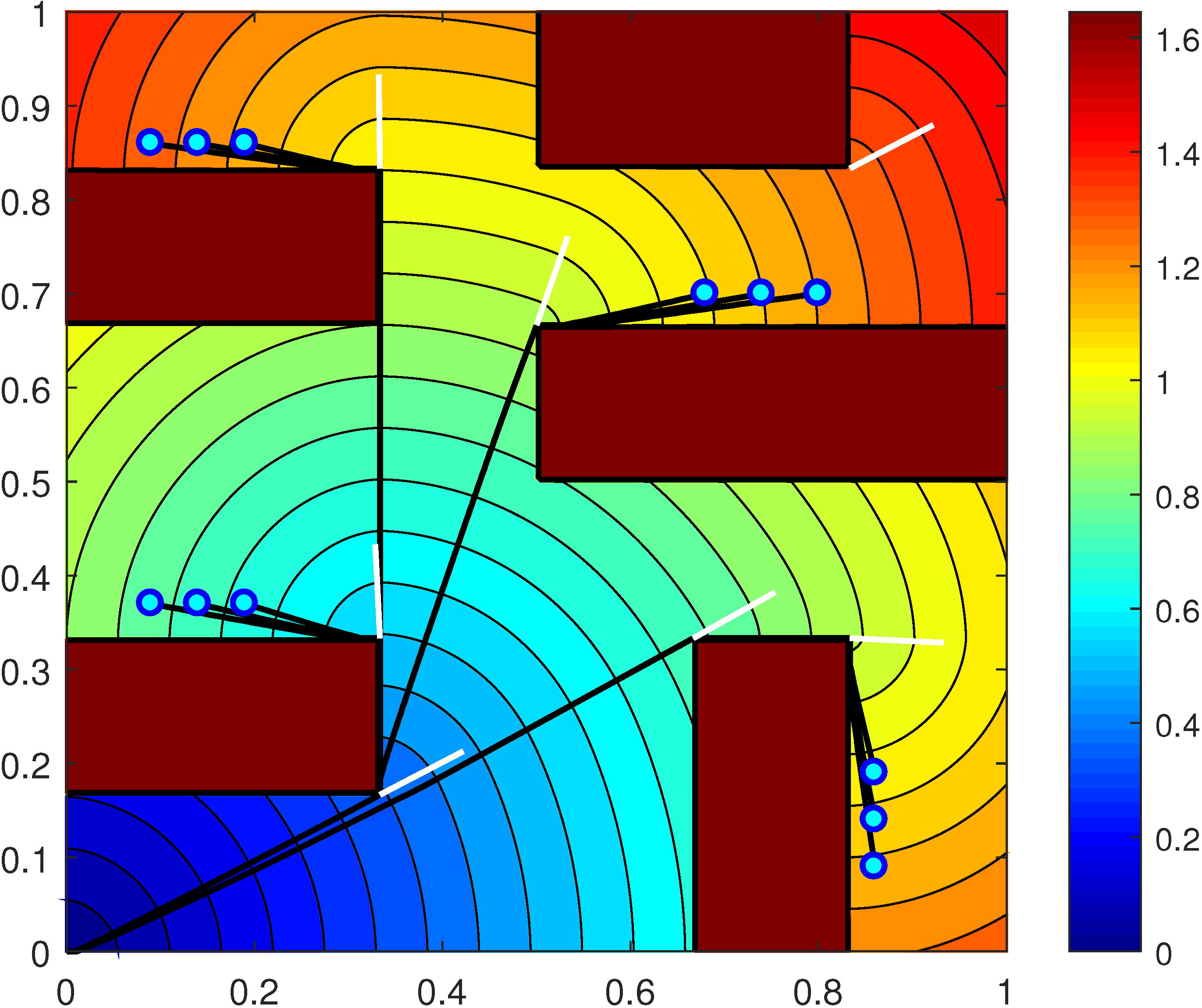}}
\subfigure[maze: localized factoring]{
\label{f:maze_1_v2}
\includegraphics[width=0.47\textwidth,clip]{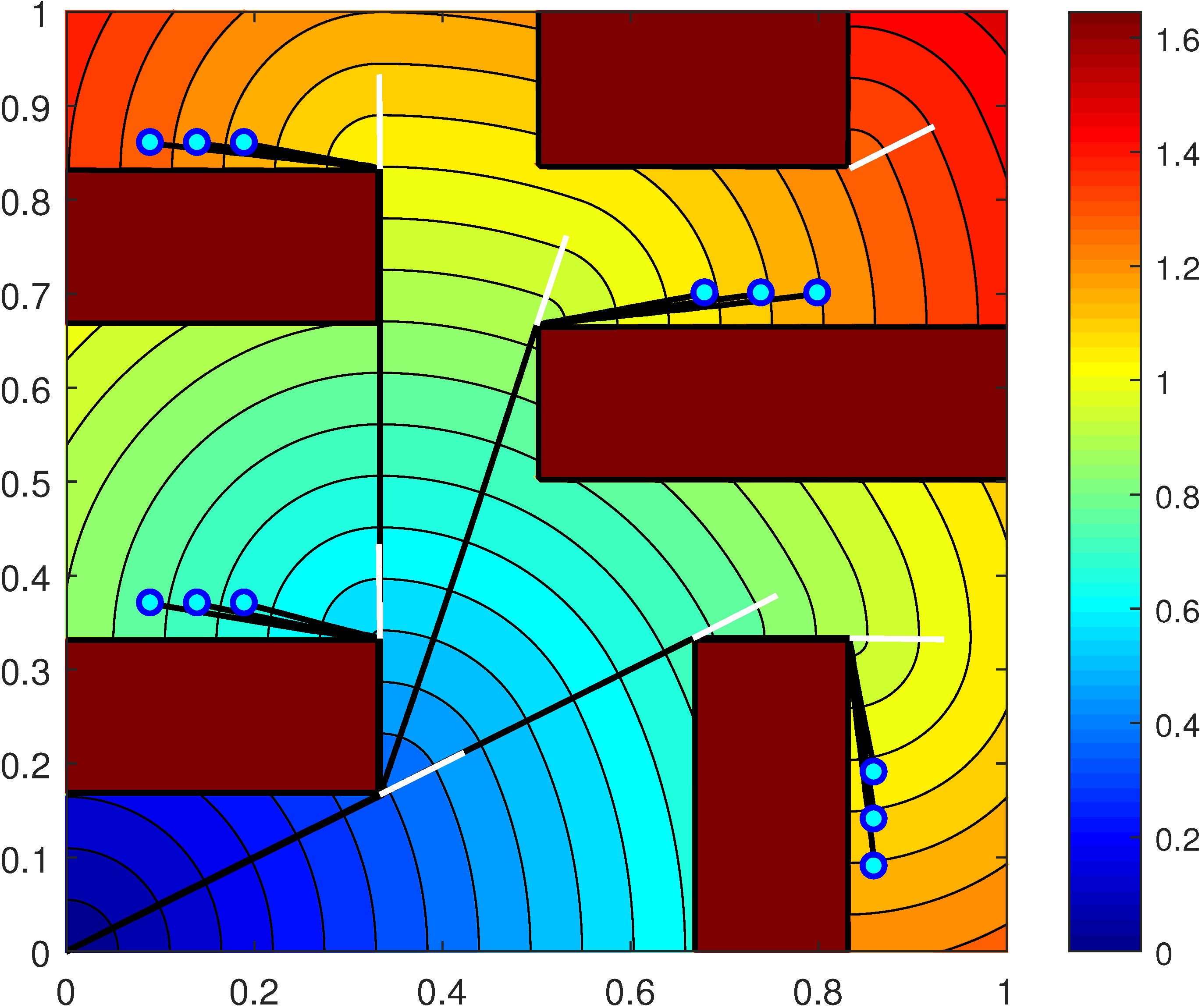}}
\caption{Navigating a maze with $F=1$.  Level sets of $u$ and representative optimal trajectories computed by the original FMM (Left)
and by Algorithm \ref{alg:dynamic} (Right).  The latter avoids obvious numerical artifacts near rarefying corners.}
\label{f:maze_1}
\end{figure}

\begin{figure}[!htb]
\centering
\subfigure[maze: original]{
\label{f:maze_2_v1}
\includegraphics[width=0.47\textwidth,clip]{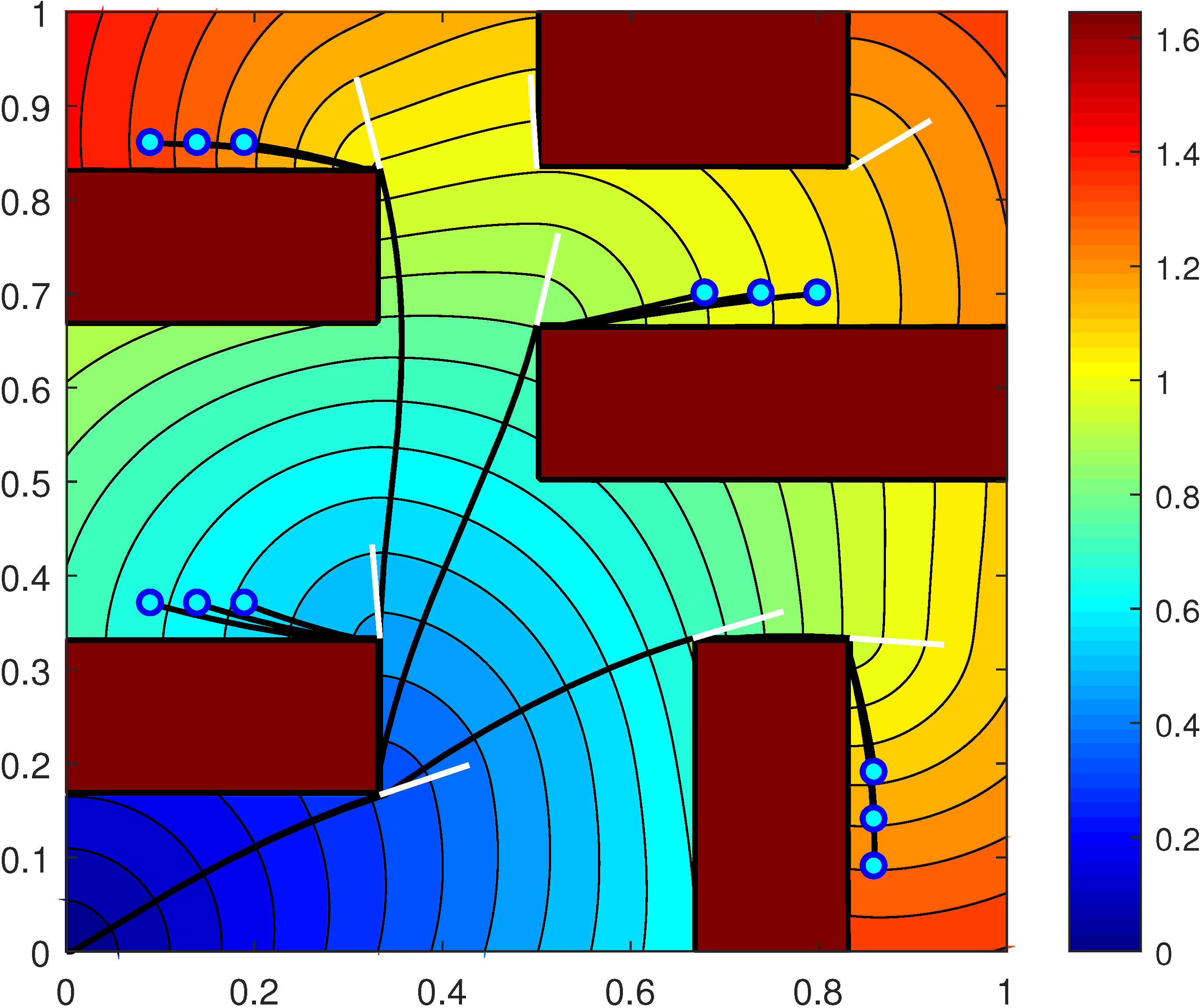}}
\subfigure[maze: localized factoring]{
\label{f:maze_2_v2}
\includegraphics[width=0.47\textwidth,clip]{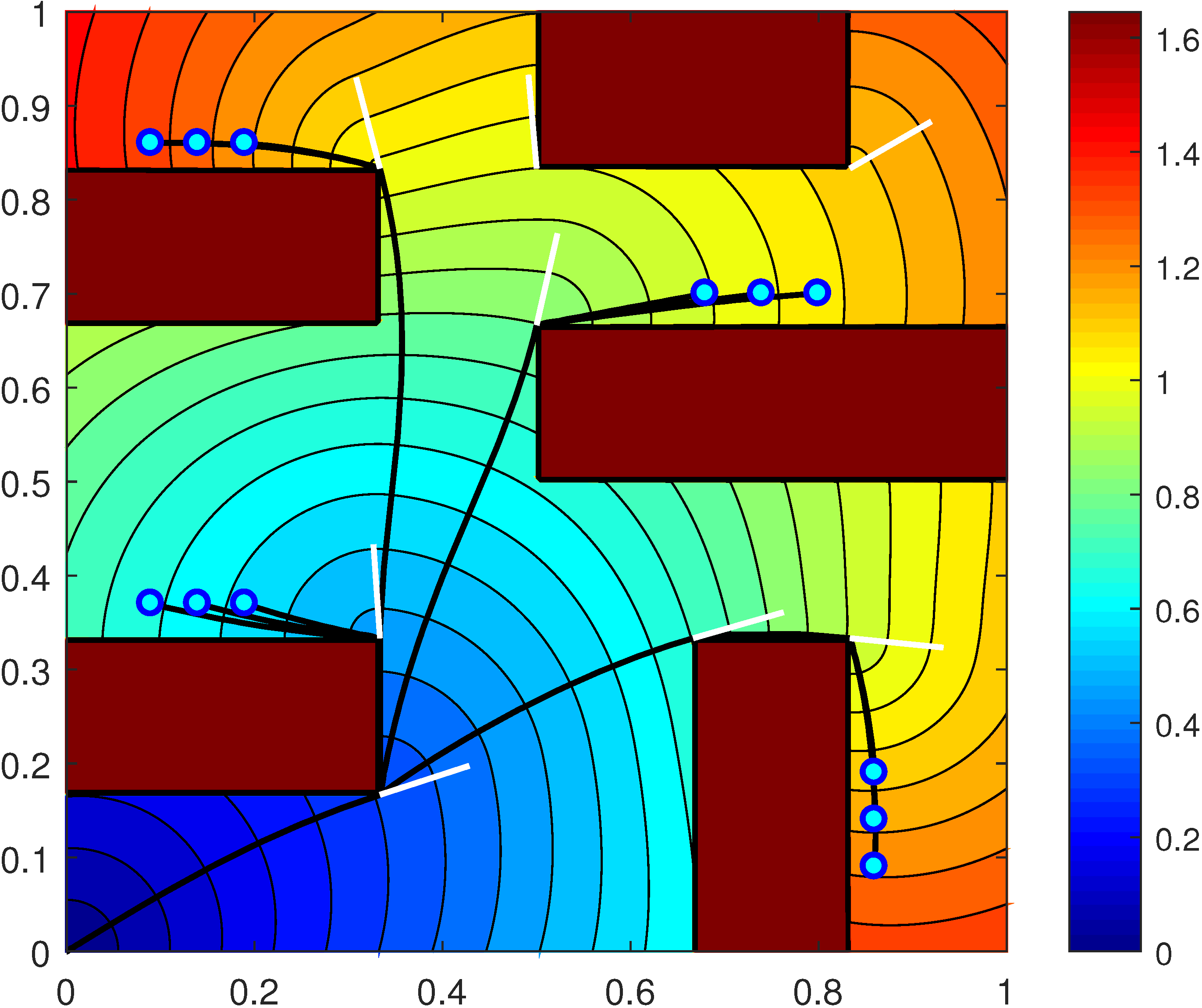}}
\caption{Navigating a maze with $F(x,y) = 1 + 0.3\sin{(2\pi x)}\sin{(2\pi y)}.$  
Level sets of $u$ and representative optimal trajectories computed by the original FMM (Left)
and by Algorithm \ref{alg:dynamic} (Right).  The latter avoids obvious numerical artifacts near rarefying corners.}
\label{f:maze_2}
\end{figure}

\begin{figure}[!htb]
\centering
\subfigure[maze with constant $F$]{
\includegraphics[width=0.28\textwidth,clip]{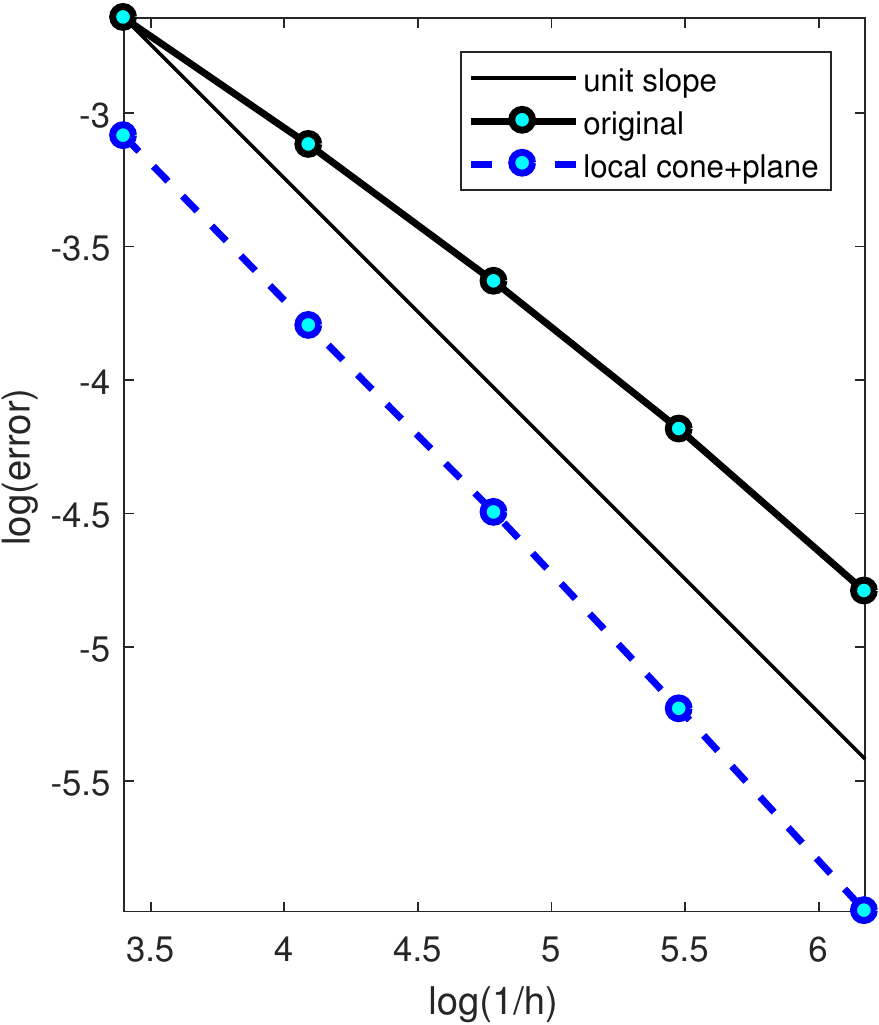}}
\hspace*{8mm}
\subfigure[maze with variable $F$]{
\includegraphics[width=0.28\textwidth,clip]{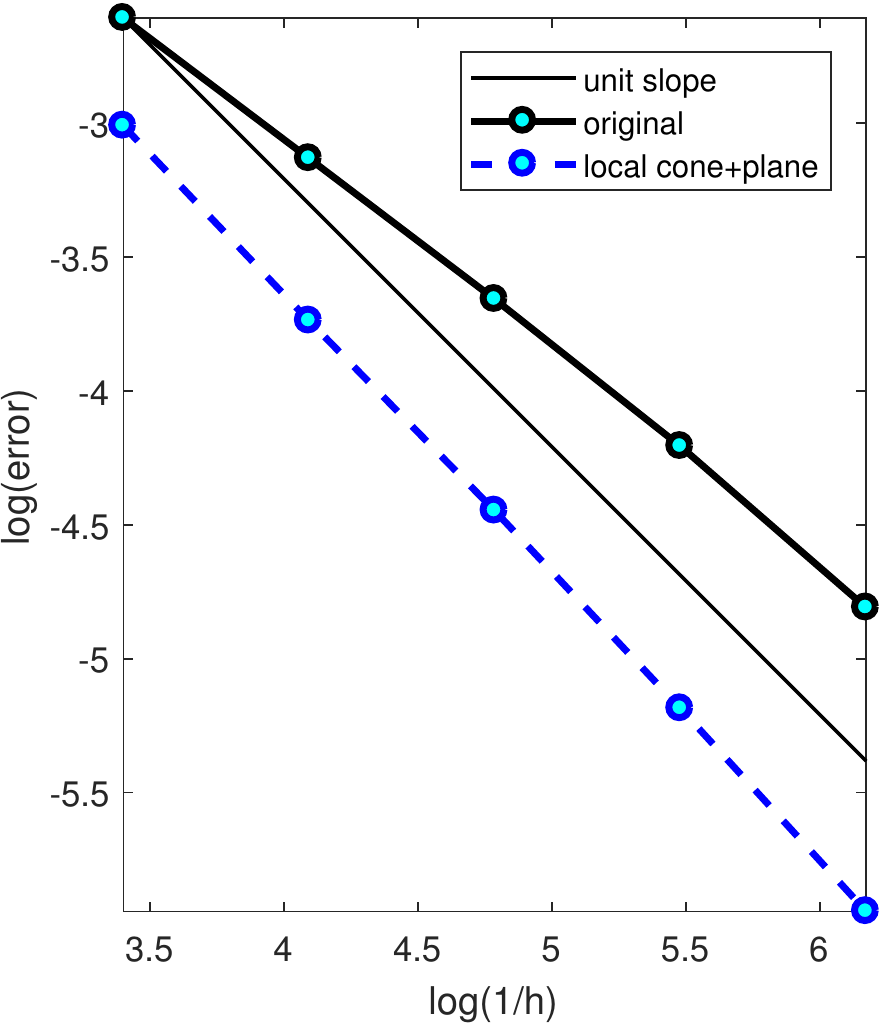}}
\caption{$L^{\infty}$ error for maze navigation examples.}
\label{f:maze_c}
\end{figure}

Since we do not rely on knowing ahead of time which corners are rarefying, 
the just-in-time localized factoring algorithm is excellent for ``maze-navigation problems'' with numerous obstacles
and possibly inhomogeneous speed function $F$. 
Before running the algorithm, we create a binary array to identify all gridpoints contained inside obstacles.
This array is then used to identify obstacle corners in Algorithm \ref{alg:dynamic}, and whenever a corner is found to be rarefying by Definition \ref{def:badcorner}, a factoring procedure is locally applied with an appropriately chosen ``additive factor'' $T(\x)$. 
Below we show two examples based on a ``maze'' from Figure \ref{f:tra}.  In Figure \ref{f:maze_1} we explore the version with $F=1.$ 
In Figure \ref{f:maze_2} we use an inhomogeneous speed $F(x,y) = 1 + 0.3\sin{(2\pi x)}\sin{(2\pi y)}.$  In both cases, the white lines are used to indicate the (local) boundaries of corner-induced rarefaction fans.  Twelve sample trajectories are shown to demonstrate that the trajectory distortions near the corners are avoided by just-in-time factoring.  The convergence is tested using the gridsizes $h = \frac{1}{30}2^{-k}$, where $k=0, \ldots, 4,$ 
and the ``ground truth'' is computed on a much finer grid with $h = 1/4800.$  Figure \ref{f:maze_c} shows that, unlike the Fast Marching Method for the original Eikonal, our approach is globally first-order accurate in both examples.

%-----------------------------------------------------------------

% !TEX root = local_factoring_for_Journal.tex

\section{Discontinuous Speed Function}
\label{ss:discont}

Rarefaction fans can also arise due to discontinuities in $F$.  
Here we consider a simple subclass of such problems:  a generalization of maze navigation examples from the previous section
to account for ``slowly permeable obstacles''.  We will assume that obstacles are described
by an open set $\Omega_{ob} \subset \left( \Omega \backslash Q \right)$ and the speed $F$ is lower inside them.
In the simplest setting, $F$ is piecewise constant with a discontinuity on $\partial   \Omega_{ob}$ and $0 < F_{ob} < F_{free}.$
We will use $\Upsilon = F_{free} / F_{ob}$ to measure the severity of obstacle slowdown. 

\begin{figure}[!htb]
\centering 
\subfigure[]{
\label{f:dis3}
\iftoggle{ForJournal}{\includegraphics[scale=.67,clip]{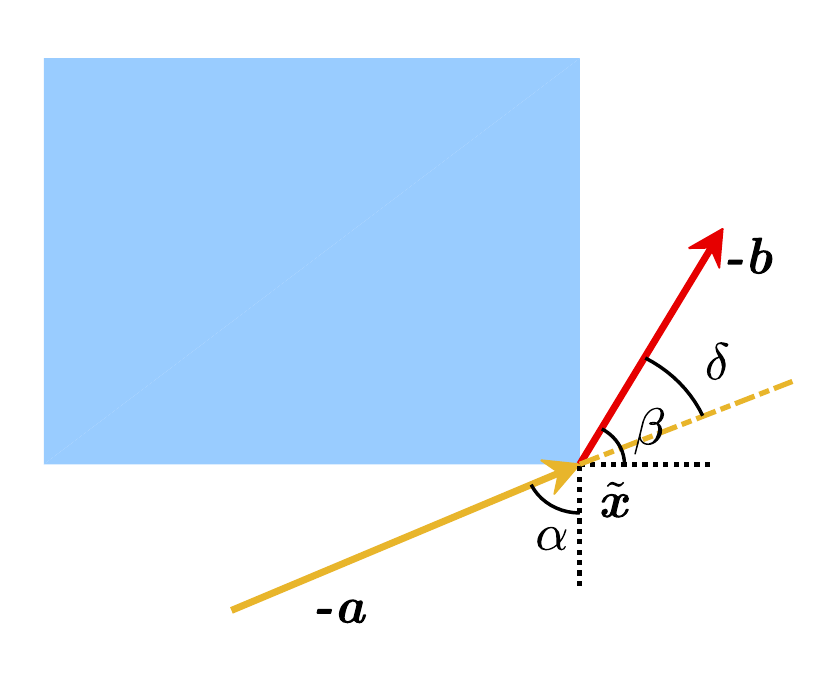}}{\includegraphics[scale=.70,clip]{discont_1}}}
\subfigure[]{
\label{f:dis4}
\iftoggle{ForJournal}{\includegraphics[scale=.67,clip]{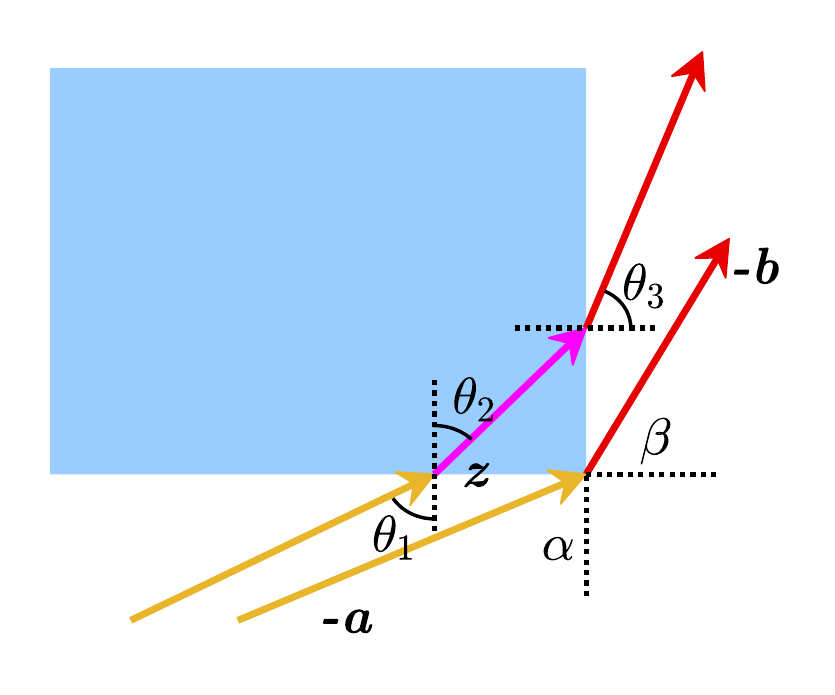}}{\includegraphics[scale=.72,clip]{discont_2}}}
\caption{A rarefaction fan at the corner of a slowly permeable obstacle. (Left) $\alpha$ is the incidence angle of a ray from the point source to the rarefying corner $\xtilde$ and $\beta$ is the ``refracted'' angle from the normal on the other side of the obstacle.
The rarefaction fan appears in a sector corresponding to the angle $\delta$ between $(-\bb)$ shown in red and the yellow dash-dotted line corresponding to $(-\ba).$ 
(Right) 
A ray refraction happening at a point $\z$ close to $\xtilde$. As $\z \to \xtilde$, $\theta_1$ and $\theta_3$ will tend to $\alpha$ and $\beta$ respectively,
with the purple segment disappearing, and yielding the $(-\ba, -\bb)$ path through $\xtilde$ in the limit.}
\label{f:dis}
\end{figure}

The following properties are relatively easy to prove for this simple type of discontinuous $F$ in 2D problems:
\begin{enumerate}
\item Rarefaction fans can only arise at point sources or at rarefying corners of slowly permeable obstacles.  
(E.g., there are no fans arising on non-corner parts of obstacle boundaries.)
\item When a rarefaction fan arises at an obstacle corner $\xtilde$, it does not propagate into that obstacle.
\item Such rarefaction fans are always confined to a sector between the characteristic direction $(-\ba)$
and another vector $(-\bb)$ found from Snell's law.
\item Suppose $\ba$ makes an angle $\alpha \in (0, \pi/2)$ with a normal to one side of a rectangular obstacle at $\xtilde$
and $(-\bb)$ makes an angle $\beta \in (0, \pi/2]$ with a normal to the other side of that obstacle; see Figure \ref{f:dis}. 
Then these angles must satisfy
\begin{equation}
\label{eq:snell}
\sin{\beta} \; = \; \min \left( \sqrt{\Upsilon^{2}-\sin^{2}{\alpha}}, \; 1 \right).
\end{equation}
and the rarefaction fan takes place in a sector of angle $\delta = (\alpha + \beta - \tfrac{\pi}{2}).$
\end{enumerate}
For the sake of brevity, we sketch the proof of the last of these only.

\begin{proof}
We can reinterpret the characteristics as light rays traveling from a point source and refracted at the boundary of slowly permeable obstacles.  We then use the Snell's Law to determine their changing directions. 

Consider a point $\z$ close to the corner $\xtilde$, whose characteristic has an incidence angle $\theta_1$ and refraction angle $\theta_2$;
see Figure \ref{f:dis4}. The incidence angle of the second refraction is $\frac{\pi}{2}-\theta_2$ and the second refraction angle is $\theta_3$. 
If $\theta_3 < \tfrac{\pi}{2}$, 
by Snell's Law these three angles must satisfy
\begin{equation}
\label{eq:angles}
\frac{\sin \theta_1 }{F_{free}} \; = \; \frac{\sin \theta_2}{F_{ob}}, 
\qquad \qquad 
\frac{\cos \theta_2 }{F_{ob}} \; = \; \frac{\sin \theta_3 }{F_{free}}.
\end{equation}
Eliminating $\theta_2$, we obtain  
\[ \sin{\theta_3} \; = \; \sqrt{\left(\frac{F_{free}}{F_{ob}}\right)^2-\sin^2{\theta_1}} \; = \; \sqrt{\Upsilon^2-\sin^2{\theta_1}} \]
We note that this equality only makes sense if $\sqrt{\Upsilon^2-\sin^2{\theta_1}} \leq 1.$
Otherwise, we will observe the ``total internal reflection'' with $\theta_3 = \tfrac{\pi}{2}$ and Snell's Law not holding for the $\theta_2$ - $\theta_3$ transition.  So, the more accurate version of  this relationship is
$\sin{\theta_3} \, = \, 
\min \left( \sqrt{\Upsilon^2-\sin^2{\theta_1}}, \, 1 \right).$
As $\textbf{z} \rightarrow \xtilde,$ we have $\theta_1 \rightarrow \alpha, \, \theta_3 \rightarrow \beta,$ and we recover \eqref{eq:snell} in the limit.
\end{proof}

\begin{rem}
It is easy to provide a sufficient condition for the rarefaction fan filling the whole region between $(-\ba)$ and the obstacle boundary (exactly as we saw in the non-permeable case).
Whenever $\Upsilon \geq \sqrt{2}$, we have $\sin \beta = 1$ and hence $\beta = \frac{\pi}{2};$ so, 
optimal trajectories from all starting positions in $\Omega \backslash\Omega_{ob}$ reach the exit set $Q$ without passing through $\Omega_{ob}.$   
On the other hand, in the continuous case ($\Upsilon = 1$), formula \eqref{eq:snell} implies that $\beta = \frac{\pi}{2} - \alpha, \delta=0$ and no rarefaction fan is present.
\end{rem}

\begin{figure}[!htb]
\centering
\subfigure[level sets of $u$]{
\label{f:perm_illustration_1}
\iftoggle{ForJournal}{\includegraphics[scale=.34,clip]{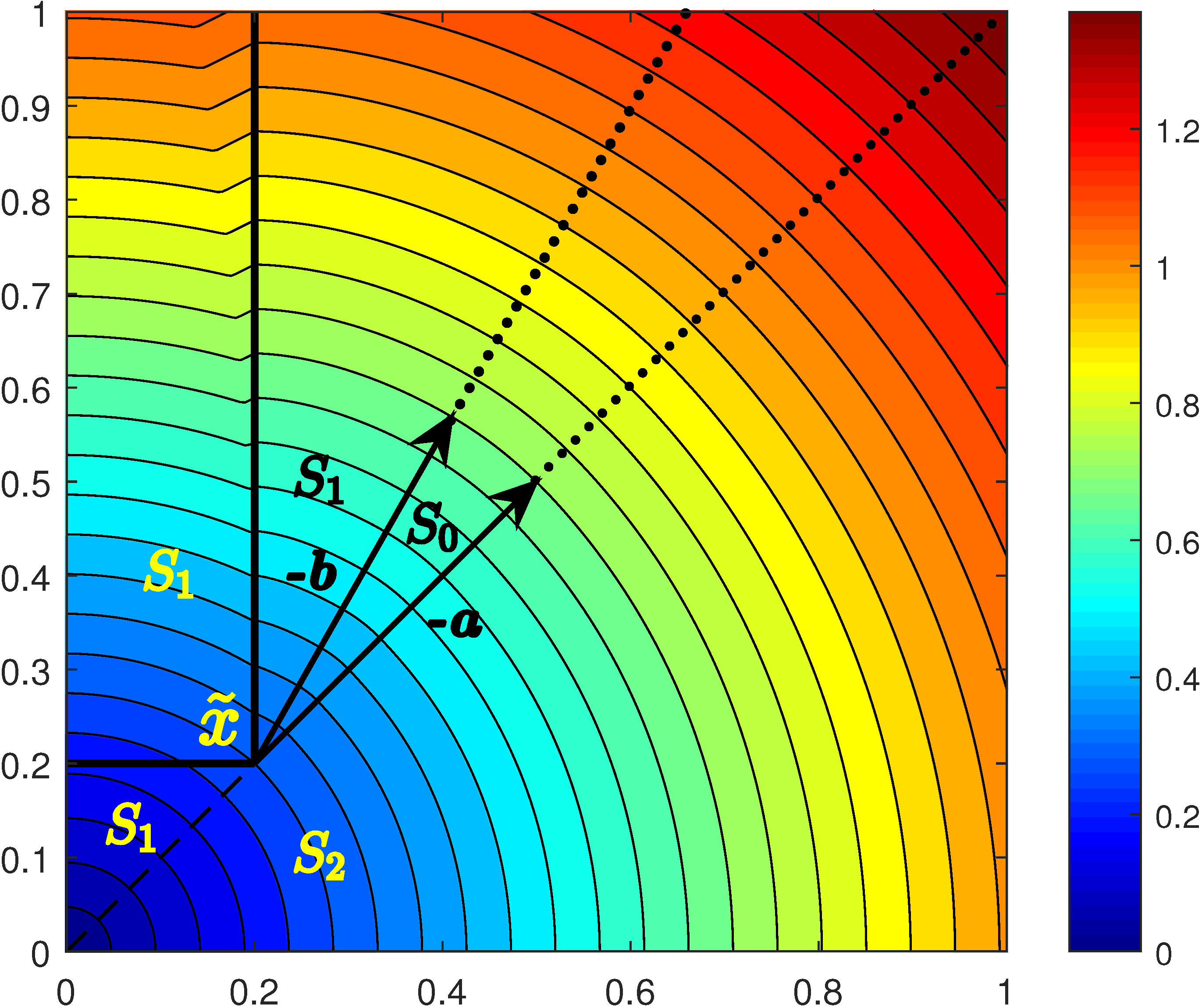}}{\includegraphics[scale=.43,clip]{perm_illustration}}}
\subfigure[domain splitting]{
\label{f:perm_illustration_2}
\iftoggle{ForJournal}{\includegraphics[scale=.34,clip]{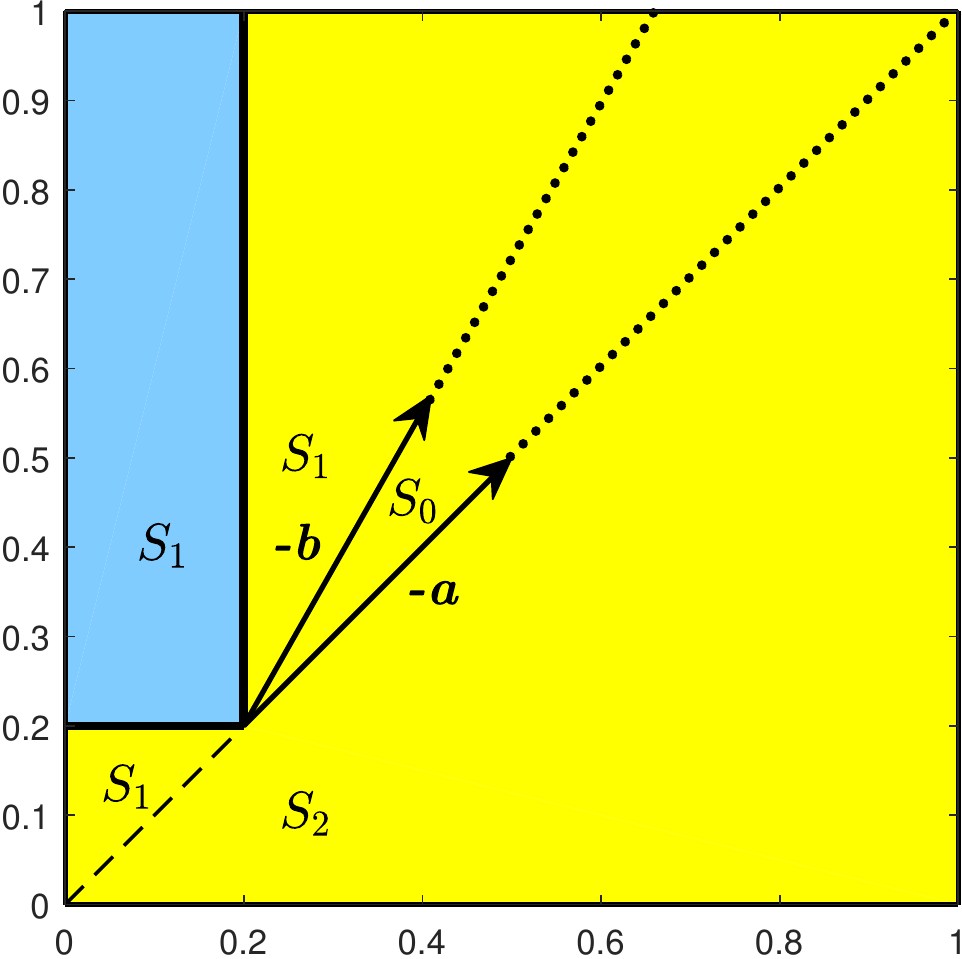}}{\includegraphics[scale=.43,clip]{perm_regions}}}
\subfigure[level sets of $T$]{
\label{f:perm_illustration_3}
\iftoggle{ForJournal}{\includegraphics[scale=.34,clip]{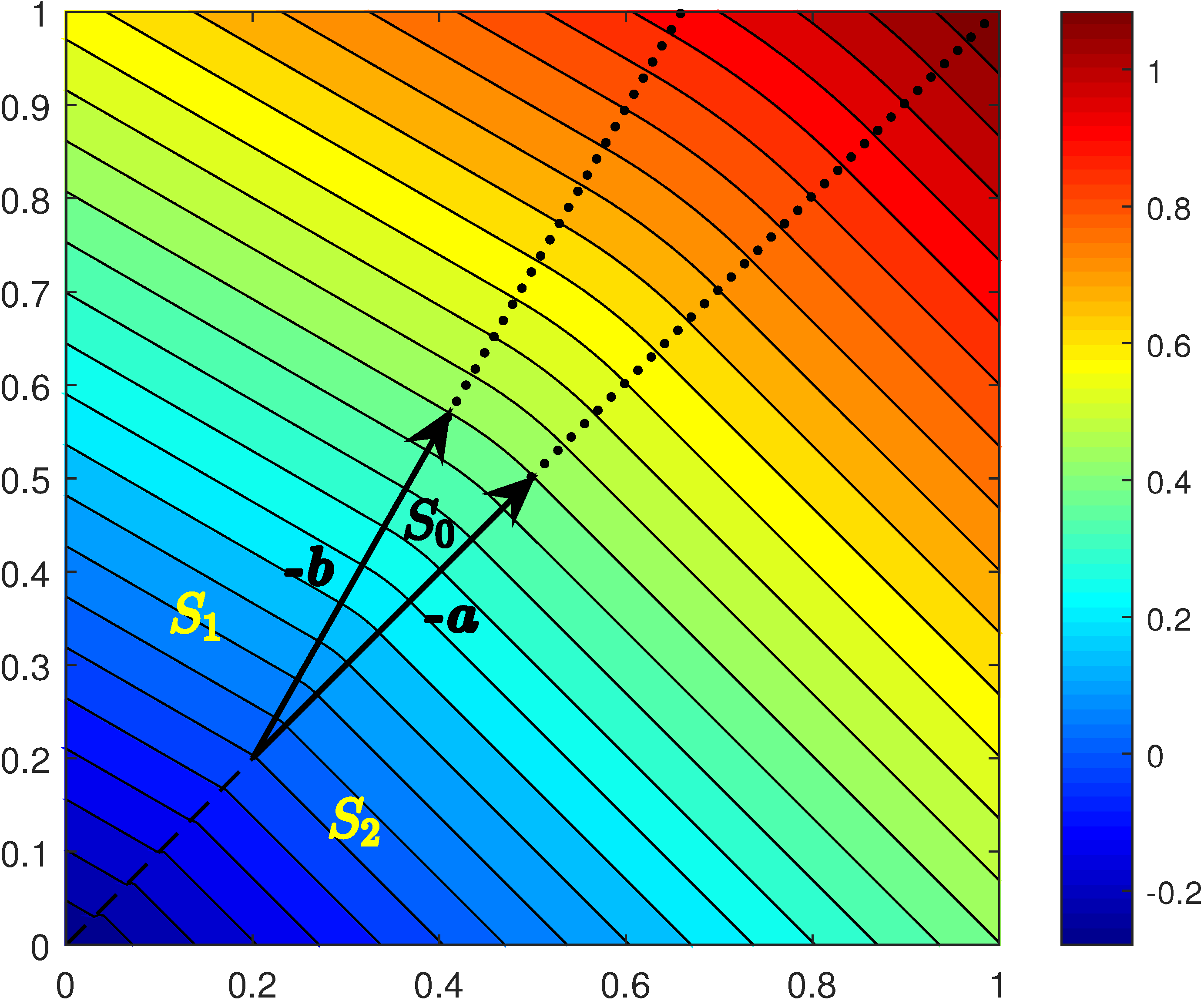}}{\includegraphics[scale=.43,clip]{perm_t_value}}}
\caption{A simple example with one ``permeable obstacle''.
(Left) level curves of the ``partially refracted'' distance to a point source.
(Center) a dynamic domain splitting based on the rarefaction fan.  
(Right) the level sets of a ``cone + 2 planes'' function $T$ capturing the correct rarefaction behavior. 
The function has a small discontinuous jump along the ray parallel to $\ba$ through $\xtilde$.}
\label{f:permeable_illustration}
\end{figure}

Using $\ba$ and $\bb$ defined at a rarefying corner $\xtilde$  in the above properties, it is natural to split $B_r(\xtilde)$ into three regions:
$$S_0 = \left\{ \x \in \Omega \mid (\x -\xtilde) \text{ is between $(-\bb)$ and $(-\ba)$ }\right\};$$
$$S_1 = \left\{ \x \in \Omega \mid (\x -\xtilde) \text{ is between $(-\bb)$ and $\ba$ }\right\};
\quad \quad
S_2 = \Omega \backslash (S_0\cup S_1).$$
We can now build a suitable (localized) factoring function $T$ as follows:
\begin{equation}
\label{eq:perm_factorT}
T(\x) \; =\;
\begin{dcases}
\dfrac{\left| \x-\xtilde \right|}{F(\xtilde)}, \quad & \x \in S_{0} \\
\dfrac{-\textbf{b} \cdot (\x-\xtilde)}{F(\xtilde)}, \quad & \x \in S_1 \\
\dfrac{-\ba \cdot (\x-\xtilde)}{F(\xtilde)}, \quad & \x \in S_2.
\end{dcases}
\end{equation}
Based on the shape of the graph, we refer to this function as a ``cone+2 planes''; see Figure \ref{f:perm_illustration_3}.
This formulation makes $T$ discontinuous along a ray parallel to $\ba = - \nabla u (\xtilde) / |\nabla u (\xtilde)|,$ 
but for a sufficiently small $r$ all gridpoints close to this ray in $B_r(\xtilde)$ will be already Accepted by the time 
we start this factoring.  Both $T$ and $\nabla T$ are continuous along $(-\ba)$ and $(-\bb).$ 

\subsection{Numerical Examples}
\label{ss:perm_example}
\begin{figure}[!htb]
\centering
\subfigure[level sets of $u(\x)$]{
\label{f:perm_1_v}
\iftoggle{ForJournal}{\includegraphics[scale=.48,clip]{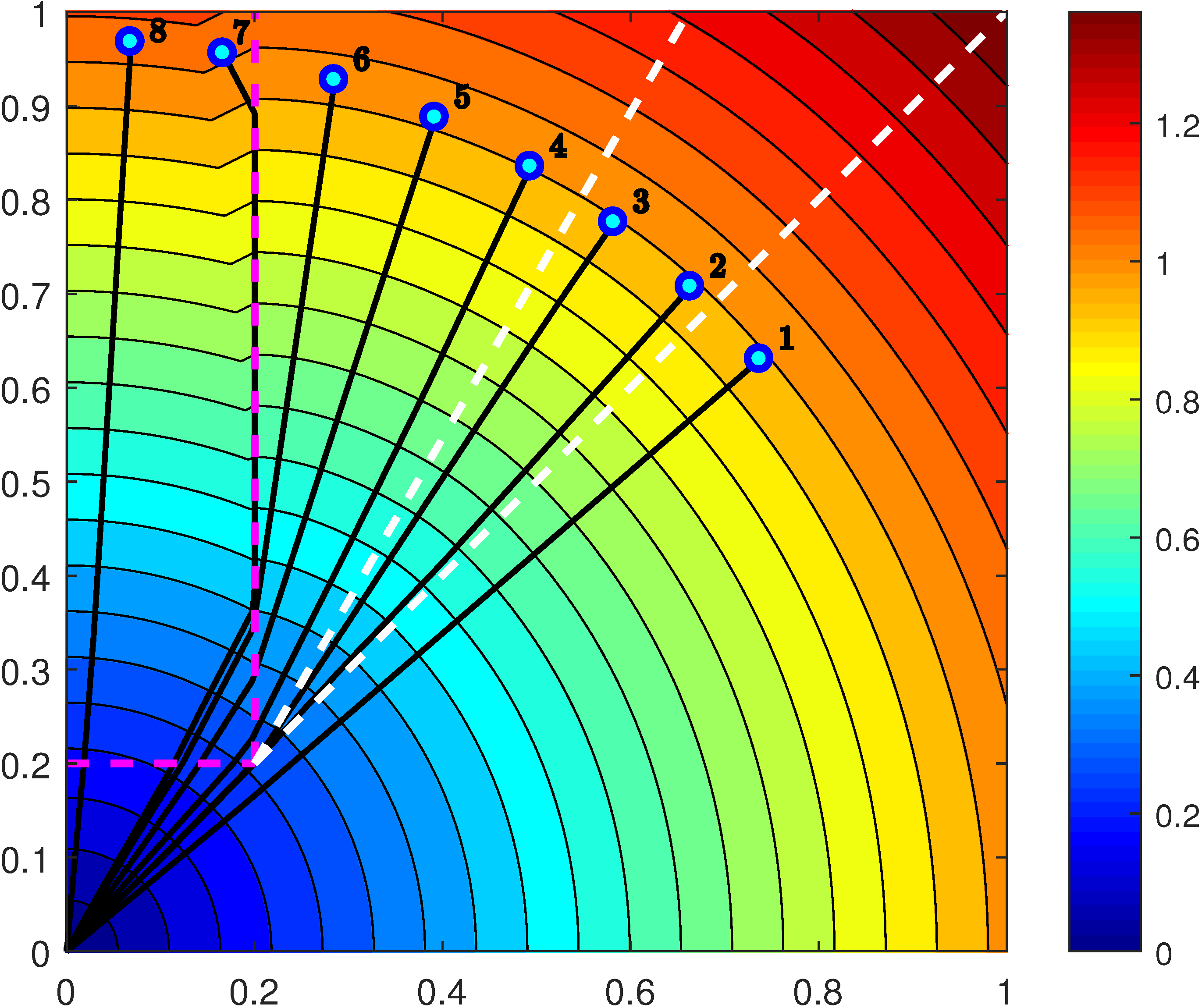}}{\includegraphics[scale=.6,clip]{perm_1_v}}}
\subfigure[$L^\infty$ error]{
\label{f:perm_1_c}
\iftoggle{ForJournal}{\hspace*{3mm} \includegraphics[scale=.35,clip]{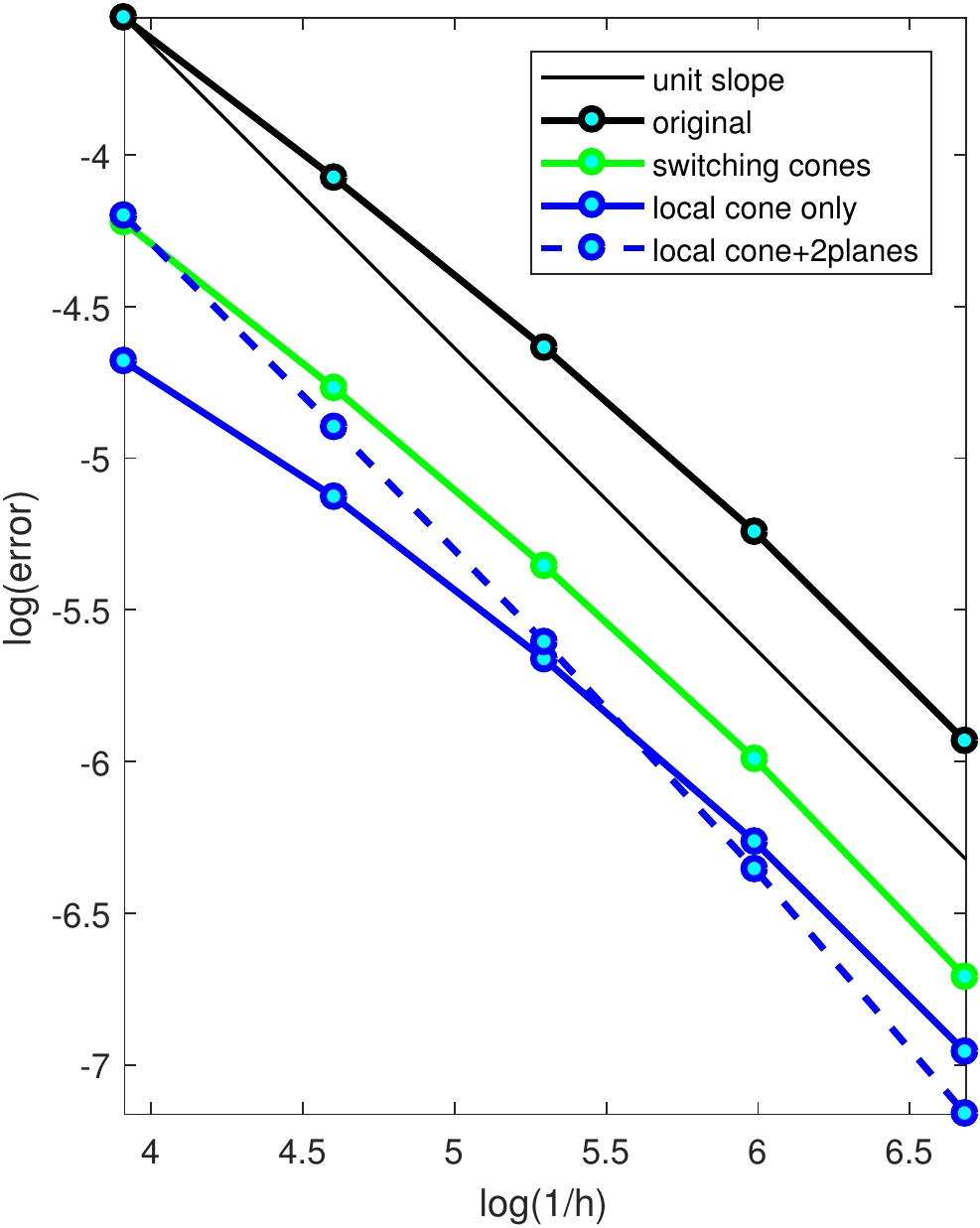}}{\includegraphics[scale=.43,clip]{perm_1_c}}}
\caption{Optimal trajectories and convergence for a ``single permeable obstacle'' example introduced in Figure \ref{f:permeable_illustration}.
(Left) The level sets of $u,$ with purple dashed lines showing the obstacle boundaries and white dashed lines showing the rarefaction fan boundaries. Eight representative optimal trajectories shown in black: 
(1) outside any region influenced by the obstacle, taking a straight line to the point source;
(2,3) starting within the rarefaction fan, coinciding after reaching the rarefying corner;
(4,5,6) experiencing a double refraction; 
(7) starting inside the obstacle and experiencing the ``total internal reflection'' described in the proof, with two different segments inside the obstacle;
(8) starting inside the obstacle and experiencing a single refraction. 
[Note: the ``light rays'' (5-8) enter the obstacle with small incidence angles, resulting in barely changed refracted angles, so the first refraction is difficult to identify visually.]}
\label{f:perm_1}
\end{figure}

Returning to the example in Figure \ref{f:permeable_illustration}, we choose $F_{ob}=\frac{2}{\sqrt{5}} \approx 0.894$ inside the obstacle 
$\Omega_{ob} =(0,0.2)\times(0.2,1)$ and $F_{free} = 1$ on $\Omega \backslash \Omega_{ob}$. 
Based on the properties discussed above, this will result in a rarefaction fan spreading in a sector of angle $\delta = \frac{\pi}{12}$ between the two white dashed lines in Figure \ref{f:perm_1_v}. 
We test the convergence of several methods described in section \ref{ss:nonpermeable_examples} 
and report the results in Figure \ref{f:perm_1_c}. 
The numerical experiments are conducted using gridsizes $h = \frac{1}{50} 2^{-k},$ where $k = 0,\ldots,4$ and the ``ground truth'' is computed on a much finer grid with $h = 1/4000$.
Unsurprisingly, the ``original'' (unfactored) method results in the largest errors and only the ``localized cone + 2 planes'' method exhibits the first-order convergence.

\begin{figure}[!htb]
\centering
\subfigure[level sets of $u(\x)$]{
\label{f:perm_2_v}
\includegraphics[width=0.45\textwidth,clip]{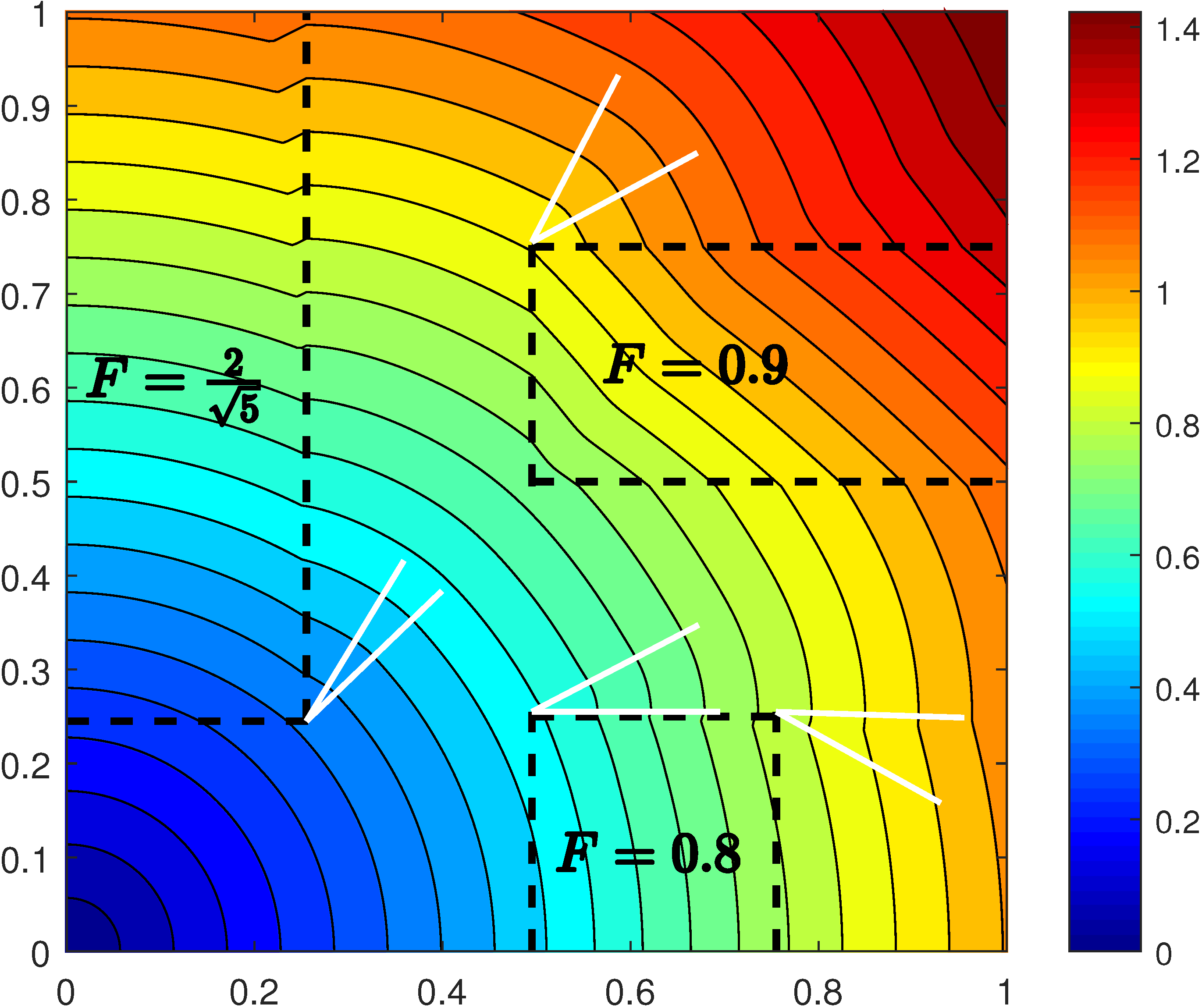}}
\subfigure[$L^\infty$ error]{
\label{f:perm_2_c}
\includegraphics[width=0.35\textwidth,clip]{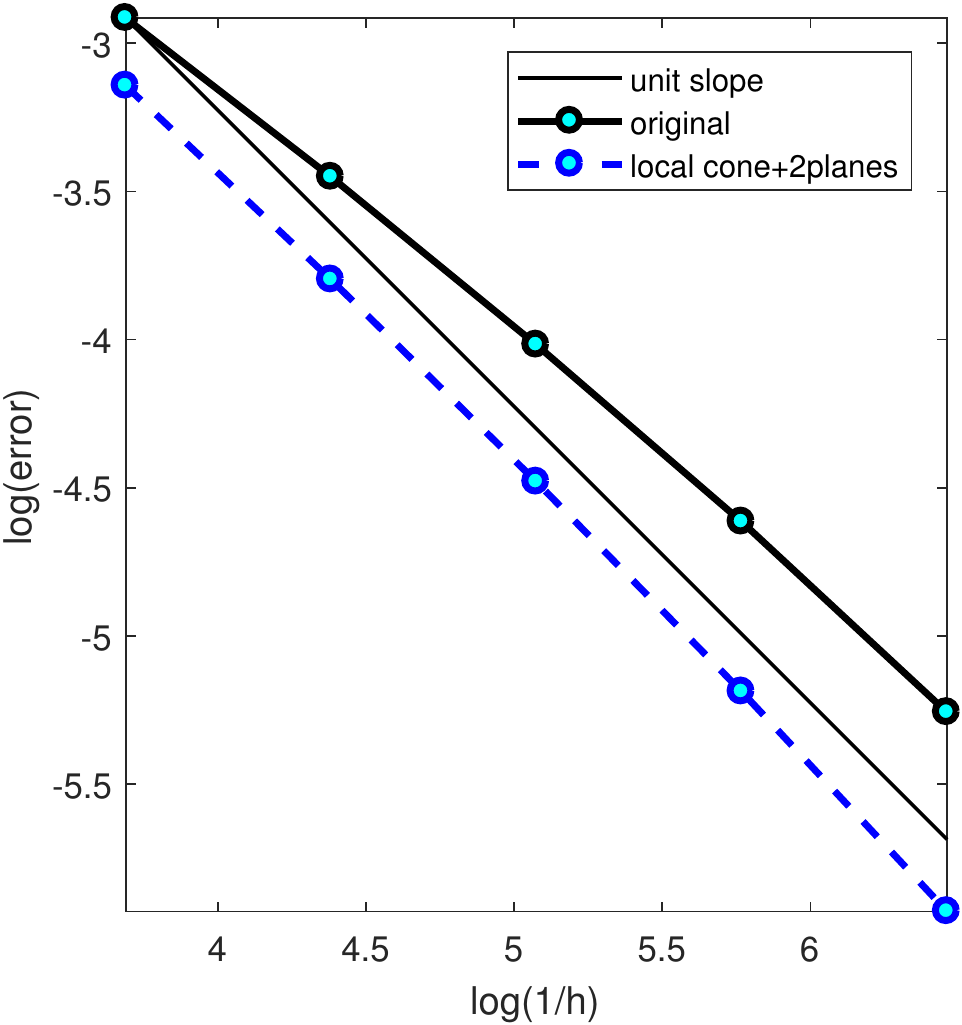}}
\caption{A maze with several slowly permeable obstacles.
(Left) The level sets of $u$ with dashed lines showing the obstacle boundaries and white line segments showing the rarefaction fan boundaries. 
The speed $F_{ob}$ is shown inside each obstacle. }
\label{f:perm_2}
\end{figure}

Our final example in Figure \ref{f:perm_2} has multiple slowly permeable obstacles with each having a different $F_{ob}$ (indicated in Figure \ref{f:perm_2_v}) and $F_{free} = 1$ in the complement.
At each corner, we use equation \eqref{eq:snell} to find $(-\bb)$ and use two white line segments to indicate the rarefying region. 
The ``ground truth'' is computed using $h = 1/6400$ and the convergence is tested using gridsizes $h = \frac{1}{40} 2^{-k}, k = 0,\ldots,4.$
Figure \ref{f:perm_2_c} demonstrates that our method reduces the errors and recovers the first-order convergence.

%-----------------------------------------------------------------

% !TEX root = local_factoring_for_Journal.tex

\section{Conclusions}
\label{s:conclude}

We have introduced a new just-in-time factoring algorithm for Eikonal equations to reduce the numerical errors due to rarefaction fans.
Prior (global and localized) factoring algorithms were meant to deal with rarefactions arising at point sources and we have 
carefully compared their accuracy in that setting.   However, our main focus has been on rarefactions arising in 2D due to nonsmothness of $\partial \Omega$ (e.g., corners of non-permeable obstacles) or discontinuities in the speed function (e.g., corners of ``slowly-permeable'' obstacles).   The locations and the geometry of such rarefaction fans are a priori unknown.  Our algorithm uncovers them dynamically and adoptively applies the localized factoring.  This dynamic aspect makes our approach natural in the Fast Marching framework. (With Fast Sweeping, one could in principle solve the original Eikonal on the entire domain, then identify all rarefaction fans in post-processing and re-solve the correctly factored equation on $\Omega.$)       
Numerical tests confirm that our method restores the full linear convergence and prevents numerical artifacts in approximating optimal trajectories once the value function is already computed.  While we have only implemented and tested the ``additive'' dynamic factoring, we expect that in the ``multiplicative'' case the results would be qualitatively similar.
All presented examples were in the context of time-optimal path planning, but other optimization criteria (e.g., a cumulative exposure to an enemy observer) would also lead to a similarly factored Eikonal equation as long as the running cost remains isotropic.
Even though our focus so far has been on applications in robotic navigation and computational geometry, we hope that the same general approach might also be useful in seismic imaging problems, which motivated much of the prior work on Eikonal factoring. 

Several straightforward generalizations will make our method more useful in practice.  
\begin{enumerate}
\item
We can easily treat general polygonal obstacles by adding dynamic factoring to prior Fast Marching techniques on (obstacle-fitted) triangulated meshes \cite{KimmSethTria, SethVlad1}.  The definition of our ``additive factor'' $T$ will stay exactly the same; see also Remark \ref{rem:discont} in section \ref{s:bad_corners}.
\item
The examples presented in section \ref{ss:discont} are based on rectangular ``slowly permeable obstacles'' with a piecewise constant speed function.  However, the same approach is also applicable for the general discontinuous speed functions as long as the discontinuity lines are polygonal and aligned with the discretization mesh.  The rarefaction fans can be determined based on a local information only (i.e., the directional limits of the speed function at a rarefying corner of the discontinuity line), and the definition of  $T$ in dynamic factoring will remain the same even when $F$ is not piecewise constant. 
\item
If the speed of motion is anisotropic (i.e., dependent on the direction of motion rather than just the current location), the value function satisfies a more general Hamilton-Jacobi-Bellman PDE.  Point-source-based factoring for the latter has already been developed (e.g., by Fast Sweeping in \cite{luo2012fast}).  Marching-type techniques for anisotropic problems (e.g., \cite{SethVlad3} or \cite{Mirebeau3}) can be similarly modified to handle the corner-induced rarefactions.     
\item
Another easy extension is to treat rarefaction fans due to more general boundary conditions (e.g., fast-varying $u=g$ specified on $\partial \Omega$ can result in rarefactions even if the boundary is smooth).  
\end{enumerate}
It will be more difficult to move to factoring suitable for higher-order accurate discretizations.
For point-source-induced fans, this has been addressed in \cite{luo2011factored} and \cite{luo2014high}.  
Similar ideas might work in our context, but higher derivatives will need to be estimated at rarefying corners and one would need to construct a smoother  $T$ than the version used in this paper.   

Finally, the obvious limitation of our current approach is that $\Omega \subset R^2.$
We expect that Eikonal problems in higher dimension will be much harder to factor dynamically.  
Even with $F=1$ and simple non-permeable box-obstacles in 3D, one would already need to deal with {\em rarefying edges} rather than corners.

%-----------------------------------------------------------------

\noindent
{\bf Acknowledgements:}  The authors are grateful to anonymous reviewers for their suggestions on improving this paper.

\bibliography{ref} 
\bibliographystyle{spmpsci}

\end{document}